\documentclass[10pt]{amsart}
\usepackage{amsmath, amssymb, latexsym}
\usepackage{mathrsfs}
\usepackage[all, knot]{xy}
\xyoption{arc}

\newtheorem{Thm}{Theorem}[section]
\newtheorem{Prop}[Thm]{Proposition}
\newtheorem{Cor}[Thm]{Corollary}
\newtheorem{Lem}[Thm]{Lemma}

\theoremstyle{definition}
\newtheorem{Def}[Thm]{Definition}
\newtheorem{Exa}[Thm]{Example}

\numberwithin{equation}{section}

\newcommand{\Z}{\mathbb{Z}}

\newcommand{\C}{\mathbb{C}}

\newcommand{\g}{\mathfrak{g}}

\newcommand{\f}{\tilde{f}}
\newcommand{\e}{\tilde{e}}

\newcommand{\gB}{\mathfrak{B}}

\newcommand{\Span}{ {\rm Span} }
\newcommand{\col}{ {\rm c} }
\newcommand{\tr}{{\rm tr}}
\newcommand{\Irr}{{\rm Irr}}
\newcommand{\Hom}{{\rm Hom}}
\newcommand{\End}{{\rm End}}
\newcommand{\Aut}{{\rm Aut}}
\newcommand{\qin}{{\rm in}}
\newcommand{\qout}{{\rm out}}
\newcommand{\HT}{{\rm ht}}
\newcommand{\wt}{{\rm wt}}
\newcommand{\cl}{{\rm cl}}
\newcommand{\Bad}[1]{B^{{\rm ad},#1}}
\newcommand{\Bone}[1]{B^{1, #1}}
\newcommand{\Ben}[1]{B^{n, #1}}
\newcommand{\im}{{\rm Im}}

\newcommand{\rowTableau}[4]
{
\fontsize{7}{7}\selectfont
\xy
(0,3.5)*{}; (0,-3.5)*{} **\dir{-};
(7,3.5)*{}; (7,-3.5)*{} **\dir{-};
(14,3.5)*{}; (14,-3.5)*{} **\dir{-};
(21,3.5)*{}; (21,-3.5)*{} **\dir{-};
(28,3.5)*{}; (28,-3.5)*{} **\dir{-};
(35,3.5)*{}; (35,-3.5)*{} **\dir{-};
(42,3.5)*{}; (42,-3.5)*{} **\dir{-};
(49,3.5)*{}; (49,-3.5)*{} **\dir{-};
(0,3.5)*{}; (49,3.5)*{} **\dir{-};
(0,-3.5)*{}; (49,-3.5)*{} **\dir{-};
(3.5,0)*{#1}; (10.5,0)*{\cdots}; (17.5,0)*{#2}; (24.5,0)*{\cdots}; (31.5,0)*{#3}; (38.5,0)*{\cdots}; (45.5,0)*{#4};
\endxy
\fontsize{10}{10}\selectfont
}

\newcommand{\colTableau}[4]
{
\fontsize{7}{7}\selectfont
\xy
(0,21)*{}; (7,21)*{} **\dir{-};
(0,14)*{}; (7,14)*{} **\dir{-};
(0,7)*{}; (7,7)*{} **\dir{-};
(0,0)*{}; (7,0)*{} **\dir{-};
(0,-7)*{}; (7,-7)*{} **\dir{-};
(0,-14)*{}; (7,-14)*{} **\dir{-};
(0,-21)*{}; (7,-21)*{} **\dir{-};
(0,21)*{}; (0,-21)*{} **\dir{-};
(7,21)*{}; (7,-21)*{} **\dir{-};
(3.5,17.5)*{#1}; (3.5,10.5)*{\vdots}; (3.5,3.5)*{#2}; (3.5,-3.5)*{#3}; (3.5,-10.5)*{\vdots}; (3.5,-17.5)*{#4};
\endxy
\fontsize{10}{10}\selectfont
}

\newcommand{\smallTableau}[2]
{
\fontsize{7}{7}\selectfont
\xy
(0,2.5)*{}; (0,-2.5)*{} **\dir{-};
(5,2.5)*{}; (5,-2.5)*{} **\dir{-};
(10,2.5)*{}; (10,-2.5)*{} **\dir{-};
(15,2.5)*{}; (15,-2.5)*{} **\dir{-};
(0,2.5)*{}; (15,2.5)*{} **\dir{-};
(0,-2.5)*{}; (15,-2.5)*{} **\dir{-};
(2.5,0)*{#1}; (7.5,0)*{\cdots}; (12.5,0)*{#2};
\endxy
\fontsize{10}{10}\selectfont
}

\newcommand{\smallTableauTheta}[4]
{
\fontsize{7}{7}\selectfont
\xy
(0,2.5)*{}; (0,-2.5)*{} **\dir{-};
(5,2.5)*{}; (5,-2.5)*{} **\dir{-};
(10,2.5)*{}; (10,-2.5)*{} **\dir{-};
(15,2.5)*{}; (15,-2.5)*{} **\dir{-};
(20,2.5)*{}; (20,-2.5)*{} **\dir{-};
(25,2.5)*{}; (25,-2.5)*{} **\dir{-};
(30,2.5)*{}; (30,-2.5)*{} **\dir{-};
(0,2.5)*{}; (30,2.5)*{} **\dir{-};
(0,-2.5)*{}; (30,-2.5)*{} **\dir{-};
(2.5,0)*{#1}; (7.5,0)*{\cdots}; (12.5,0)*{#2}; (17.5,0)*{#3};(22.5,0)*{\cdots};(27.5,0)*{#4};
\endxy
\fontsize{10}{10}\selectfont
}

\newcommand{\sBox}[1]
{
\fontsize{7}{7}\selectfont
\xy
(0,2.5)*{}; (0,-2.5)*{} **\dir{-};
(5,2.5)*{}; (5,-2.5)*{} **\dir{-};
(0,2.5)*{}; (5,2.5)*{} **\dir{-};
(0,-2.5)*{}; (5,-2.5)*{} **\dir{-};
(2.5,0)*{#1};
\endxy
\fontsize{10}{10}\selectfont
}

\newcommand{\ssBox}[1]
{
\fontsize{7}{7}\selectfont
\xy
(0,2.5)*{}; (0,-2.5)*{} **\dir{-};
(7,2.5)*{}; (7,-2.5)*{} **\dir{-};
(0,2.5)*{}; (7,2.5)*{} **\dir{-};
(0,-2.5)*{}; (7,-2.5)*{} **\dir{-};
(3.5,0)*{#1};
\endxy
\fontsize{10}{10}\selectfont
}

\newcommand{\Patterns}
{
\fontsize{7}{7}\selectfont
\xy
(0,0)*{}; (56,0)*{} **\dir{-};
(0,7)*{}; (56,7)*{} **\dir{-};
(0,14)*{}; (56,14)*{} **\dir{-};
(0,21)*{}; (56,21)*{} **\dir{-};
(4,0)*{}; (4,28)*{} **\dir{-};
(11,0)*{}; (11,28)*{} **\dir{-};
(18,0)*{}; (18,28)*{} **\dir{-};
(25,0)*{}; (25,28)*{} **\dir{-};
(32,0)*{}; (32,28)*{} **\dir{-};
(42,0)*{}; (42,28)*{} **\dir{-};
(49,0)*{}; (49,28)*{} **\dir{-};
(56,0)*{}; (56,28)*{} **\dir{-};
(7.5, 3.5)*{n}; (14.5, 3.5)*{0}; (21.5, 3.5)*{1}; (28.5, 3.5)*{2}; (37, 3.5)*{\cdots}; (45.5, 3.5)*{k-1}; (52.5, 3.5)*{k};
(7.5, 10.5)*{0}; (14.5, 10.5)*{1}; (21.5, 10.5)*{2}; (28.5, 10.5)*{3}; (37, 10.5)*{\cdots}; (45.5, 10.5)*{k}; (52.5, 10.5)*{k+1};
(7.5, 17.5)*{1}; (14.5, 17.5)*{2}; (21.5, 17.5)*{3}; (28.5, 17.5)*{4}; (37, 17.5)*{\cdots}; (45.5, 17.5)*{k+1}; (52.5, 17.5)*{k+2};
(7.5, 25)*{\vdots}; (14.5, 25)*{\vdots}; (21.5, 25)*{\vdots}; (28.5, 25)*{\vdots}; (45.5, 25)*{\vdots}; (52.5, 25)*{\vdots};
(30,-3)*{\text{the pattern $\mathtt{P}_k^1$ of weight $\Lambda_k$}};
(70,0)*{}; (126,0)*{} **\dir{-};
(70,7)*{}; (126,7)*{} **\dir{-};
(70,14)*{}; (126,14)*{} **\dir{-};
(70,21)*{}; (126,21)*{} **\dir{-};
(74,0)*{}; (74,28)*{} **\dir{-};
(81,0)*{}; (81,28)*{} **\dir{-};
(88,0)*{}; (88,28)*{} **\dir{-};
(95,0)*{}; (95,28)*{} **\dir{-};
(102,0)*{}; (102,28)*{} **\dir{-};
(112,0)*{}; (112,28)*{} **\dir{-};
(119,0)*{}; (119,28)*{} **\dir{-};
(126,0)*{}; (126,28)*{} **\dir{-};
(77.5, 3.5)*{2}; (84.5, 3.5)*{1}; (91.5, 3.5)*{0}; (98.5, 3.5)*{n}; (107, 3.5)*{\cdots}; (115.5, 3.5)*{k+1}; (122.5, 3.5)*{k};
(77.5, 10.5)*{1}; (84.5, 10.5)*{0}; (91.5, 10.5)*{n}; (98.5, 10.5)*{n-1}; (107, 10.5)*{\cdots}; (115.5, 10.5)*{k}; (122.5, 10.5)*{k-1};
(77.5, 17.5)*{0}; (84.5, 17.5)*{n}; (91.5, 17.5)*{n-1}; (98.5, 17.5)*{n-2}; (107, 17.5)*{\cdots}; (115.5, 17.5)*{k-1}; (122.5, 17.5)*{k-2};
(77.5, 25)*{\vdots}; (84.5, 25)*{\vdots}; (91.5, 25)*{\vdots}; (98.5, 25)*{\vdots}; (115.5, 25)*{\vdots}; (122.5, 25)*{\vdots};
(100,-3)*{\text{the pattern $\mathtt{P}_k^n$ of weight $\Lambda_k$}};
\endxy
\fontsize{10}{10}\selectfont
}

\newcommand{\Quiver}
{
\fontsize{7}{7}\selectfont
\xy
(15,8)*{\bullet}; (0,0)*{} **\dir{-} ?>* \dir{>};
(0,0)*{\bullet}; (10,0)*{\bullet} **\dir{-} ?>* \dir{>};
(20,0)*{\bullet}; (30,0)*{} **\dir{-} ?>* \dir{>};
(30,0)*{\bullet}; (15,8)*{} **\dir{-} ?>* \dir{>};
(15,10)*{0}; (0,-2)*{1}; (10,-2)*{2}; (20,-2)*{n-1}; (30,-2)*{n};
(-5,6)*{(I,\Omega):=}; (33,0)*{,};
(50,0)*{} ; (65,8)*{\bullet} **\dir{-} ?>* \dir{>} ;
(60,0)*{\bullet};(50,0)*{\bullet}  **\dir{-} ?>* \dir{>};
(80,0)*{} ; (70,0)*{\bullet}  **\dir{-} ?>* \dir{>};
(65,8)*{}; (80,0)*{\bullet}  **\dir{-} ?>* \dir{>};
(65,10)*{0}; (50,-2)*{1}; (60,-2)*{2}; (70,-2)*{n-1}; (80,-2)*{n};
(45,6)*{(I,\overline{\Omega}):=};
\endxy
\fontsize{10}{10}\selectfont
}

\newcommand{\ExYWone}
{
\fontsize{7}{7}\selectfont
\xy
(0,-5)*{}; (0,-10)*{} **\dir{-};
(5,-5)*{}; (5,-10)*{} **\dir{-};
(10,0)*{}; (10,-10)*{} **\dir{-};
(15,0)*{}; (15,-10)*{} **\dir{-};
(0,-10)*{}; (15,-10)*{} **\dir{-};
(0,-5)*{}; (15,-5)*{} **\dir{-};
(10,0)*{}; (15,0)*{} **\dir{-};
(2.5,-7.5)*{1}; (7.5,-7.5)*{2}; (12.5,-7.5)*{0}; (12.5,-2.5)*{1}; (17,-10)*{,} ;
(20,-5)*{}; (20,-10)*{} **\dir{-};
(25,-5)*{}; (25,-10)*{} **\dir{-};
(30,5)*{}; (30,-10)*{} **\dir{-};
(35,5)*{}; (35,-10)*{} **\dir{-};
(20,-10)*{}; (35,-10)*{} **\dir{-};
(20,-5)*{}; (35,-5)*{} **\dir{-};
(30,0)*{}; (35,0)*{} **\dir{-};
(30,5)*{}; (35,5)*{} **\dir{-};
(22.5,-7.5)*{1}; (27.5,-7.5)*{2}; (32.5,-7.5)*{0}; (32.5,-2.5)*{1}; (32.5, 2.5)*{2}; (37,-10)*{,} ;
(40,0)*{}; (40,-5)*{} **\dir{-};
(45,0)*{}; (45,-5)*{} **\dir{-};
(50,10)*{}; (50,-5)*{} **\dir{-};
(55,10)*{}; (55,-5)*{} **\dir{-};
(60,10)*{}; (60,-5)*{} **\dir{-};
(40,-5)*{}; (60,-5)*{} **\dir{-};
(40,0)*{}; (60,0)*{} **\dir{-};
(50,5)*{}; (60,5)*{} **\dir{-};
(50,10)*{}; (60,10)*{} **\dir{-};
(42.5,-2.5)*{1}; (47.5,-2.5)*{2}; (52.5,-2.5)*{0}; (57.5,-2.5)*{1}; (52.5,2.5)*{1}; (57.5,2.5)*{2}; (52.5,7.5)*{2}; (57.5,7.5)*{0};
\endxy
\fontsize{10}{10}\selectfont
}

\newcommand{\ExYWen}
{
\fontsize{7}{7}\selectfont
\xy
(0,0)*{}; (0,-5)*{} **\dir{-};
(5,0)*{}; (5,-5)*{} **\dir{-};
(10,10)*{}; (10,-5)*{} **\dir{-};
(15,10)*{}; (15,-5)*{} **\dir{-};
(0,-5)*{}; (15,-5)*{} **\dir{-};
(0,0)*{}; (15,0)*{} **\dir{-};
(10,5)*{}; (15,5)*{} **\dir{-};
(10,10)*{}; (15,10)*{} **\dir{-};
(2.5,-2.5)*{2}; (7.5,-2.5)*{1}; (12.5,-2.5)*{0}; (12.5,2.5)*{2}; (12.5,7.5)*{1}; (17,-5)*{,} ;
(20,0)*{}; (20,-5)*{} **\dir{-};
(25,10)*{}; (25,-5)*{} **\dir{-};
(30,10)*{}; (30,-5)*{} **\dir{-};
(20,-5)*{}; (30,-5)*{} **\dir{-};
(20,0)*{}; (30,0)*{} **\dir{-};
(25,5)*{}; (30,5)*{} **\dir{-};
(25,10)*{}; (30,10)*{} **\dir{-};
(22.5,-2.5)*{1}; (27.5,-2.5)*{0}; (27.5,2.5)*{2}; (27.5,7.5)*{1}; (32,-5)*{,} ;
(35,-5)*{}; (35,-10)*{} **\dir{-};
(40,-5)*{}; (40,-10)*{} **\dir{-};
(45,-5)*{}; (45,-10)*{} **\dir{-};
(50,0)*{}; (50,-10)*{} **\dir{-};
(55,5)*{}; (55,-10)*{} **\dir{-};
(60,5)*{}; (60,-10)*{} **\dir{-};
(35,-10)*{}; (60,-10)*{} **\dir{-};
(35,-5)*{}; (60,-5)*{} **\dir{-};
(50,0)*{}; (60,0)*{} **\dir{-};
(55,5)*{}; (60,5)*{} **\dir{-};
(37.5,-7.5)*{2}; (42.5,-7.5)*{1}; (47.5,-7.5)*{0}; (52.5,-7.5)*{2}; (57.5,-7.5)*{1};
(52.5,-2.5)*{1}; (57.5,-2.5)*{0}; (57.5,2.5)*{2};
\endxy
\fontsize{10}{10}\selectfont
}

\begin{document}

\title[Quiver varieties and adjoint crystals of level $\ell$ for type $A_n^{(1)}$]
{Quiver varieties and adjoint crystals of level $\ell$ for type $A_n^{(1)}$}
\author[Euiyong Park]{Euiyong Park $^1$ } 
\thanks{$^1$ This work was supported by KRF Grant \# 2007-341-C00001.}
\address{School of Mathematics, Korea Institute for Advanced Study,
87 Hoegiro, Dongdaemun-gu, Seoul 130-722, Korea}
\email{eypark@kias.re.kr}

\subjclass[2000]{17B10, 17B67, 81R10}
\keywords{quiver varieties, adjoint crystals, crystal graphs}

\begin{abstract}

Let $B(\Lambda)$ be a level $\ell$ highest weight crystal of the quantum affine algebra $U_q(A_n^{(1)})$.
We construct an explicit crystal isomorphism between the geometric realization $\gB(\Lambda)$ of the crystal $B(\Lambda)$
using quiver varieties and
the path realization $\mathcal{P}^{\rm ad}(\Lambda)$ arising from the adjoint crystal $\Bad{\ell}$ of level $\ell$.

\end{abstract}

\maketitle


\vskip 2em

\section*{Introduction}

The {\em theory of perfect crystals} was developed in studying vertex models
in terms of representation theory of quantum affine algebras \cite{KKMMNN91}.
The theory of perfect crystals has a lot of applications to several research areas, and plays an
important role in the theory of crystal bases. The crystal $B(\Lambda)$ of an irreducible highest weight
module $V(\Lambda)$ with highest weight $\Lambda$ over a quantum affine algebra can be realized as the crystal $\mathcal{P}^B(\Lambda)$ consisting
of $\Lambda$-paths in a perfect crystal $B$. In \cite{BFKL06}, Benkart, Frenkel, Kang and Lee gave a uniform construction
of level 1 perfect crystals $\Bad{1} = B(0) + B(\theta)$, called the {\em adjoint crystals}, for all quantum affine algebras. Here, $\theta$ is the 
maximal (short) root of the corresponding finite-dimensional simple Lie algebras. It was expected that the crystal
$ \Bad{\ell} = B(0) + B(\theta)+ \cdots + B(\ell \theta) $ is a perfect crystal of level $\ell$.
This conjecture was proved for types $A_n^{(1)}, C_n^{(1)}, A_{2n}^{(2)}, D_{n+1}^{(2)}$ and $D_{n}^{(1)}$ (which yield types $B_n^{(1)}, A_{2n-1}^{(2)}$)
 in \cite{KKM94,KKMMNN92,Kde09,AS06}.
These perfect crystals are called the {\em adjoint crystals of level} $\ell$.

Let $\g$ be a symmetric Kac-Moody algebra. Lusztig \cite{L90, L91} gave a geometric construction of $U_q^-(\g)$ in terms of perverse sheaves and introduced
the notion of {\em canonical bases} (or {\em lower global bases} by Kashiwara \cite{Kash91}). In \cite{KS97}, Kashiwara and Saito proved that the set $\gB(\infty)$ of irreducible components
of Lusztig's quiver varieties has a crystal structure, which is isomorphic to the crystal $B(\infty)$ of $U_q^-(\g)$.
For a dominant integral weight $\Lambda$, Nakajima \cite{N94, N98} defined a new family of quiver varieties corresponding to $\Lambda$ and constructed
a geometric realization of $V(\Lambda)$. In \cite{St02}, Saito defined a crystal structure on the set $\gB(\Lambda)$ of irreducible
components of certain Lagrangian subvarieties of Nakajima's quiver varieties, and showed that $\gB(\Lambda)$ is isomorphic
to the crystal $B(\Lambda)$.

In this paper, we give an explicit crystal isomorphism from the geometric realization $\gB(\Lambda)$ to the path realization $\mathcal{P}^{\rm ad}(\Lambda)$
arising from the adjoint crystal $\Bad{\ell}$ of level $\ell$ for type $A_n^{(1)}$. This is a generalization of the previous result \cite{KP10}.
Let $\g$ be the Kac-Moody algebra of type $A_n^{(1)}$ and $\overset{\circ}{\g}$ be the corresponding simple Lie algebra.
Let $\Bone{\ell}$ (resp.\ $\Ben{\ell}$) be the perfect crystal of level $\ell$ given by $\eqref{Eq: def of B1l}$ (resp.\ $\eqref{Eq: def of Bnl}$).
Choose the counterclockwise orientation $\Omega$ of the Dynkin diagram of type $A_n^{(1)}$. Let $\Lambda$ be a dominant integral weight of level $\ell$, and let
$\Lambda' := \sum_{i\in I} \Lambda(h_{i+1})\Lambda_i$ and $\Lambda'' := \sum_{i\in I} \Lambda(h_{i-1})\Lambda_i$, where $\Lambda_i$ is the $i$th fundamental weight.
We first deduce two kinds of {\it Young walls} \cite{HK02,Kang03} from the results in \cite{FLOTW99, JMMO91, KL06},
and describe irreducible components of quiver varieties by using these two kinds of Young walls;
see Proposition \ref{Prop:isom between Yn and gB}, Proposition \ref{Prop:isom between Y1 and gB}.
Then, using these descriptions,
we give a geometric interpretation of the {\em fundamental isomorphism of perfect crystals}
\begin{align*}
\gB(\Lambda) \simeq \gB(\Lambda')\otimes \Bone{\ell}\ \text{ and }\ \gB(\Lambda) \simeq \gB(\Lambda'')\otimes \Ben{\ell}
\end{align*}
(Theorem \ref{Thm:Fundamantal thm for qv}).
Unlike the approach in \cite{KP10}, we deal directly with irreducible components in the crystal $\gB(\Lambda)$ without using the strict embedding $\gB(\Lambda) \hookrightarrow \gB(\infty) \otimes T_\Lambda \otimes C$.
From these fundamental isomorphisms, we give explicit isomorphisms from $\gB(\Lambda)$ to the
path realizations $\mathcal{P}^1(\Lambda)$ and $\mathcal{P}^n(\Lambda)$ arising from $\Bone{\ell}$ and $\Ben{\ell}$ in terms of dimension vectors (Theorem \ref{Thm:iso from gB to P1 and Pn}).
Using the crystal isomorphism $ \Bad{\ell} \simeq \Bone{\ell} \otimes \Ben{\ell} $, we obtain a geometric interpretation of the isomorphism
\begin{align*}
\gB(\Lambda) \simeq \gB(\Lambda)\otimes \Bad{\ell}
\end{align*}
and an explicit map
$$ \gB(\Lambda) \buildrel \sim\over\longrightarrow  \mathcal{P}^{\rm ad}(\Lambda) $$
in terms of dimension vectors; see
Theorem \ref{Thm:iso from gB to Pad} and Corollary \ref{Cor: fundamental thm for ad}.

This paper is organized as follows. Section 1 contains the notion of crystal for $U_q(A_n^{(1)})$ and the theory of perfect crystals.
Section 2 gives a description of the adjoint crystals $\Bad{\ell}$ of level $\ell$ for type $A_n^{(1)}$ in terms of {\em Young tableaux}.
In Section 3, we review the geometric realizations $\gB(\infty)$ and $\gB(\Lambda)$ using quiver varieties.
In Section 4, we deduce two kinds of Young walls $\mathcal{Y}^1(\Lambda)$ and $\mathcal{Y}^n(\Lambda)$ which give combinatorial realizations of $B(\Lambda)$,
and investigate connections between $\mathcal{Y}^1(\Lambda)$ (resp.\ $\mathcal{Y}^n(\Lambda)$) and $\gB(\Lambda)$.
Then we describe irreducible components of quiver varieties by using these two kinds of Young walls.
In Section 5, we construct an explicit crystal isomorphism $\gB(\Lambda) \buildrel \sim\over\longrightarrow  \mathcal{P}^{\rm ad}(\Lambda)$.
We first find geometric interpretations of the isomorphisms
$\gB(\Lambda) \simeq \gB(\Lambda')\otimes \Bone{\ell}$ and $\gB(\Lambda) \simeq \gB(\Lambda'')\otimes \Ben{\ell} $
from the description using Young walls in Section 4.
Then we give an explicit description of the $\Lambda$-path in $\Bone{\ell}$ (resp.\ $\Ben{\ell}$) corresponding to an irreducible component $X \in \gB(\Lambda)$ in
terms of dimension vectors. Using the crystal isomorphism $ \Bad{\ell} \simeq \Bone{\ell} \otimes \Ben{\ell} $, we obtain
an explicit isomorphism $\gB(\Lambda) \buildrel \sim \over \longrightarrow \mathcal{P}^{\rm ad}(\Lambda)$ in terms of dimension vectors.

\vskip 1em

\noindent
{\bf Acknowledgments.} The author would like to express his deepest gratitude to Professor Seok-Jin Kang for his guidance, valuable comments and enthusiastic encouragement.

\vskip 2em

\section{Crystal graphs for $U_q(A_n^{(1)})$} \label{Sec:crystals}

Let $I:=\Z/(n+1)\Z$ be an index set. The {\it affine Cartan
datum of type $A_n^{(1)}$} consists of (i) the {\em affine Cartan matrix} $A = (a_{ij})_{i,j\in I}$ of type $A_n^{(1)}$,
(ii) {\em dual weight lattice} $P^{\vee} := \bigoplus_{i=0}^{n}\Z h_i \oplus
\Z d$, (iii) {\em affine weight lattice} $P := \bigoplus_{i=0}^n \Z \Lambda_i
\oplus \Z \delta \subset \mathfrak{h}^*$, where
$\mathfrak{h} := \C \otimes_\Z P^{\vee},\ \Lambda_i(h_j) = \delta_{ij},\ \Lambda_i(d)=0,\
\delta(h_i)=0,\ \delta(d)=1\ (i,j\in I),$
(iv) the set of {\em simple coroots} $\Pi^{\vee} = \{ h_i|\ i\in I \}$,
(v) the set of {\em simple roots} $\Pi = \{ \alpha_i|\ i\in I \}$ given by
$$\alpha_j(h_i) = a_{ij}\ \text{ and } \alpha_j(d)=\delta_{0,j}\ (i,j\in
I).$$
The free abelian group $Q := \bigoplus_{i=0}^n \Z \alpha_i$
is called the {\it root lattice} and the semigroup $Q^+ :=
\sum_{i=0}^n \Z_{\ge0} \alpha_i$ is called the {\it positive root
lattice}. For $\alpha = \sum_{i\in I} k_i \alpha_i \in Q^+$,
the number $\HT(\alpha) := \sum_{i\in I} k_i$ is called the {\it height} of $\alpha$.
The elements in $ P^+ := \{ \lambda \in P \mid \lambda(h_i) \ge 0,\ i\in I
\} $ are called the {\it dominant integral weights}.

Let $U_q(\g)$ be the {\it quantum affine algebra} associated with the affine Cartan
datum $(A, \Pi, \Pi^{\vee}, P, P^{\vee})$ of type $A_n^{(1)}$.
Let $U^-_q(\g)$ (resp.\ $U^+_q(\g)$) be the subalgebra of $U_q(\g)$ generated by $f_i$ (resp.\ $e_i$) for $ i\in I$.

\begin{Def} A {\it crystal} $B$ associated with $U_q(\g)$ is a set
together with the maps $\wt: B \to P,\ \tilde{e}_i, \tilde{f}_i: B \to
B \sqcup \{0\},$ and $\varepsilon_i, \varphi_i: B \to \Z \cup \{ -\infty \}\ (i\in I)$ satisfying
the following conditions:
\begin{enumerate}
\item $\varphi_i(b) = \varepsilon_i(b) + \langle h_i, \wt(b) \rangle$ for all $i\in I$,
\item $\wt(\tilde e_i b) = \wt(b)+\alpha_i$ if $\tilde{e}_ib \in B$,
\item $\wt(\tilde f_i b) = \wt(b)-\alpha_i$ if $\tilde{f}_ib \in B$,
\item $\varepsilon_i(\tilde e_i b)= \varepsilon_i(b)-1,\ \varphi_i(\tilde e_i b) = \varphi_i(b)+1$ if $\tilde e_i b \in B$,
\item $\varepsilon_i(\tilde f_i b)= \varepsilon_i(b)+1,\ \varphi_i(\tilde f_i b) = \varphi_i(b)-1$ if $\tilde f_i b \in B$,
\item $\tilde f_i b = b'$ if and only if $b=\tilde e_i b'$ for $b,b'\in B,\ i\in I$,
\item if $\varphi_i(b)=- \infty$ for $b \in B$, then $\tilde e_i b = \tilde f_i b = 0$.
\end{enumerate}
\end{Def}
We refer the reader to \cite{HK02,Kash90} for the {\it tensor product rule} of crystals.
For $\Lambda \in P^+$, let $V(\Lambda)$ be the irreducible highest weight $U_q(\g)$-module with highest weight $\Lambda$.
It was shown in \cite{Kash91} that $U_q^-(\g)$ and $V(\Lambda)$ have crystal bases, denoted by $B(\infty)$ and $B(\Lambda)$ respectively.

We now recall the notion of perfect crystals. Let $U'_q(\g)$ be the subalgebra of $U_q(\g)$ generated by $e_i,
f_i, q^{\pm h_i}\ (i\in I)$, and set
$\overline{P} := \bigoplus_{i=0}^n\Z \Lambda_i$,
$\overline{P}^+ := \sum_{i=0}^n\Z_{\ge 0} \Lambda_i$,
$\overline{P}^{\vee} :=
\bigoplus_{i=0}^n\Z h_i$ and $\overline{\mathfrak{h}} := \C \otimes_\Z
\overline{P}^{\vee}$.
Denote by $\cl:P \to \overline{P}$ the natural projection from $P$ to $\overline{P}$.
Given a $U'_q(\g)$-crystal $B$ and $b \in B$, let
$$ \varepsilon(b) := \sum_{i=0}^n \varepsilon_i(b)\Lambda_i,
\qquad \varphi(b) :=\sum_{i=0}^n \varphi_i(b)\Lambda_i. $$

\begin{Def} A {\it perfect crystal of level $\ell$} is a finite
$U_q'(\g)$-crystal $B$ satisfying the following conditions:
\begin{enumerate}
\item there exists a finite dimensional $U'_q(\g)$-module with a crystal basis whose crystal graph is isomorphic to $B$,
\item $B \otimes B$ is connected,
\item there exists a classical weight $\Lambda' \in \overline{P}$ such that
$$ \wt(B) \subset \Lambda' + \sum_{i \ne 0} \Z_{\le 0} \alpha_i, \qquad \#(B_{\Lambda'} := \{ b\in B \mid \wt(b) = \Lambda' \})=1,$$
\item for any $b\in B$, we have $\langle c:=\sum_{i\in I} h_i,\ \varepsilon(b) \rangle \ge \ell$,
\item for each $\Lambda \in \overline{P}_{\ell}^+:= \{ \mu \in \overline{P}^+|\ \langle c,\
\mu \rangle = \ell \}$, there exist unique vectors $b^{\Lambda}$
and $b_{\Lambda} $ in $B$ such that $\varepsilon(b^\Lambda) =
\Lambda,\ \varphi(b_\Lambda) = \Lambda$.
\end{enumerate}
\end{Def}

Given $\Lambda \in P^+$ with level $\ell$
and a perfect crystal $B$ of level $\ell$, it was shown in
\cite{KKMMNN91, KKMMNN92} that there exists a unique crystal isomorphism,
called the {\it fundamental isomorphism of perfect
crystals},
\begin{align} \label{Eq:fundamental thm of perfect crystals}
\Phi: B(\Lambda) \overset{\sim} \longrightarrow
B(\varepsilon(b_{\Lambda})) \otimes B
\end{align}
sending the highest weight vector $u_{\Lambda}$ to
$u_{\varepsilon(b_{\Lambda})} \otimes b_{\Lambda}$. By applying this
crystal isomorphism repeatedly, we get a sequence of crystal
isomorphisms
$$B(\Lambda) \overset{\sim} \longrightarrow B(\Lambda^1) \otimes
B \overset{\sim} \longrightarrow B(\Lambda^2) \otimes B \otimes B
\overset {\sim} \longrightarrow \cdots\cdots,$$ where
$\Lambda^0=\Lambda$, $b_0=b_{\Lambda}$,
$\Lambda^{k+1}=\varepsilon(b_k)$, $b_{k+1}=b_{\Lambda^{k+1}}$
$(k\ge 0)$. The sequence $\mathbf{p}_{\Lambda} :=
(b_k)_{k=0}^\infty$ is called the {\it ground-state path of weight
$\Lambda$} and a sequence $\mathbf{p}=(p_k)_{k=0}^\infty$ of
elements $p_k \in B$ is called a {\it $\Lambda$-path in $B$} if
$p_k = b_k$ for all $k \gg 0$. Let ${\mathcal P}^{B}(\Lambda)$ be the set of $\Lambda$-paths in $B$.
Then we have the {\it path realization} of $B(\Lambda)$ arising from $B$.

\begin{Thm}\cite{KKMMNN91} \label{Thm: crystal iso of path realization}
There exists a unique crystal isomorphism $B(\Lambda) \overset{\sim}
\longrightarrow \mathcal{P}^{ B}(\Lambda)$ which maps
$u_\Lambda$ to $\mathbf{p}_\Lambda$.
\end{Thm}

\vskip 2em

\section{Adjoint crystals of level $\ell$ for type $A_n^{(1)}$}

In this section, we give a description of adjoint crystals $\Bad{\ell}$ for type $A_n^{(1)}$ in terms of {\em Young tableaux} \cite{Kde09}.
Let $U_q(\overset{\circ}{\g})$ be the subalgebra of $U_q(\g)$ generated by $f_i, e_i$ and $ q^{h_i}$ for $i\in I \setminus \{0 \}$, and
let $\varpi_i := \Lambda_i - \Lambda_0\ (1 \le i \le n)$. It is well-known that a $U_q(\overset{\circ}{\g})$-crystal
can be realized as the set of semistandard tableaux with entries $1, \ldots, n+1$ (\cite[Chapter 7]{HK02}, \cite{KN94}, etc).


Let us first recall the perfect crystals $\Bone{\ell}$ and $\Ben{\ell}$ given in \cite{KKMMNN92}.
Note that, though $\Bone{\ell}$ and $\Ben{\ell}$ are dual to each other, it is difficult to describe an explicit crystal isomorphism between the path realizations arising from $\Bone{\ell}$ and $\Ben{\ell}$.
 As a $U_q(\overset{\circ}{\g})$-crystal,
$\Bone{\ell}$ is isomorphic to the crystal $B(\ell \varpi_1)$. For $b \in B(\ell \varpi_1)$, let $\nu_i(b)$ be the number of $i$ in $b$, and set
$\nu_{0}(b) := \nu_{n+1}(b)$. Then
\begin{align} \label{Eq: def of B1l}
\Bone{\ell} = \left\{ b :=\rowTableau{1}{1}{n+1}{n+1}
 \  \mid\  \sum_{i=1}^{n+1} \nu_i(b) = \ell \right \}.
\end{align}
The crystal structure of $\Bone{\ell}$ is given in \cite[Section 1.2]{KKMMNN92}.
The map $\wt : \Bone{\ell} \to \wt(\Bone{\ell})$ is bijective. More precisely, let
\begin{align} \label{Eq:map wt to B1}
\psi^{1,\ell}: \wt(\Bone{\ell}) \longrightarrow \Bone{\ell}
\end{align}
be the map sending a weight $\sum_{k\in I} a_k\Lambda_k $ to the tableau $b$ given by
$$  \nu_i(b) = \frac{1}{n+1} \left(\ell - \sum_{k=1}^n k a_k \right)  + \sum_{k=i}^{n}a_k \quad \text{ for }i = 1, \ldots, n+1. $$
Then $\psi^{1,\ell}$ is the inverse of $\wt : \Bone{\ell} \to \wt(\Bone{\ell})$. For each $\Lambda = \sum_{k\in I} a_k\Lambda_k \in P^+$ of level $\ell$,
the ground-state path $\mathbf{p}^1(\Lambda)$ of $\Lambda$ is given as follows:
\begin{align} \label{Eq:g-s path in B1}
\mathbf{p}^1(\Lambda) = ( b^{\Lambda}_k )_{k \ge 0},
\end{align}
where $b^{\Lambda}_k$ is the tableau in $\Bone{\ell}$ obtained by $\nu_j(b^{\Lambda}_k) = a_{j+k \mod n}$  for $j=1, \ldots, n+1 $.
We denote by $\mathcal{P}^1(\Lambda)$ the path realization arising from $\Bone{\ell}$.

We now consider the crystal $\Ben{\ell}$. As a $U_q(\overset{\circ}{\g})$-crystal,
$\Ben{\ell}$ is isomorphic to the crystal $B(\ell \varpi_n)$. Let
$$ \begin{tabular}{|c|}
     \hline
     $\overline{i}$ \\
     \hline
   \end{tabular}
\ :=\  \colTableau{1}{i-1}{i+1}{n+1} \ .
 $$
For $\overline{b} \in B(\ell \varpi_n)$, let $\overline{\nu}_i(\overline{b})$ denote the number of $\overline{i}$ in $\overline{b}$, and set
$\overline{\nu}_0(\overline{b}) := \overline{\nu}_{n+1}(\overline{b})$. Then
\begin{align} \label{Eq: def of Bnl}
 \Ben{\ell} = \left\{ \overline{b} := \rowTableau{\overline{n+1}}{\overline{n+1}}{\overline{1}}{\overline{1}}
 \  \mid\  \sum_{i=1}^{n+1} \overline{\nu}_i(\overline{b}) = \ell \right \}.
\end{align}
The crystal structure of $\Ben{\ell}$ is given in \cite[Section 1.2]{KKMMNN92}.
Then the map $\wt : \Ben{\ell} \to \wt(\Ben{\ell})$ is bijective, and the inverse is given by
\begin{align} \label{Eq:map wt to Bn}
\psi^{n,\ell}: \wt(\Ben{\ell}) \longrightarrow \Ben{\ell},
\end{align}
which maps a weight $\sum_{k\in I} a_k\Lambda_k $ to the tableau $\overline{b}$ obtained by
$$  \overline{\nu}_i(\overline{b}) = \frac{1}{n+1} \left(\ell + \sum_{k=1}^n k a_k \right)  - \sum_{k=i}^{n}a_k \quad \text{ for }i = 1, \ldots, n+1. $$
For $\Lambda = \sum_{k\in I} a_k\Lambda_k \in P^+ $ of level $\ell$, the ground-state path $\mathbf{p}^n(\Lambda)$ of $\Lambda$ is given as follows:
\begin{align} \label{Eq:g-s path in Bn}
\mathbf{p}^n(\Lambda) = ( \overline{b}^{\Lambda}_k )_{k \ge 0},
\end{align}
where $\overline{b}^{\Lambda}_k$ is the tableau in $\Ben{\ell}$ obtained by $\overline{\nu}_j(\overline{b}^{\Lambda}_k) = a_{j-k-1 \mod n}$  for $j=1, \ldots, n+1 $.
We denote by $\mathcal{P}^n(\Lambda)$ the path realization arising from $\Ben{\ell}$.

Let $\theta := \varpi_1 + \varpi_n$. Then any tableau $T$ of $B(k\theta)\ (1 \le k \le \ell)$ can be written as
\begin{align} \label{Eq:tableau in B(k theta)}
T = \smallTableauTheta{\overline{i}_1}{\overline{i}_k}{j_1}{j_k}
\end{align}
for some $i_1 \ge \cdots \ge i_k $ and $j_1 \le \cdots \le j_k $. As a $U_q(\overset{\circ}{\g})$-crystal, let
$$ \Bad{\ell} := \bigoplus_{k=0}^\ell B(k \theta). $$
To define the Kashiwara operators $\e_0, \f_0$ on $\Bad{\ell}$, we need some notation.
For a tableau $T$ and $ s \in \Z_{\ge0}$, let $\mathfrak{R}_s(T)$ (resp.\ $\mathfrak{L}_s(T)$) be the tableau obtained from $T$ by removing
$s$ columns from the left of $T$ (resp.\ the right of $T$). In particular, if we take a tableau $T$ of the form $\eqref{Eq:tableau in B(k theta)}$,
then we regard $\mathfrak{R}_k(T) = \smallTableau{j_1}{j_k}$ (resp.\ $\mathfrak{L}_k(T) = \smallTableau{\overline{i}_1}{\overline{i}_k}$)
as an element in the perfect crystal $\Bone{k}$ (resp.\ $\Ben{k}$).
Given two tableaux $T \in B(k \varpi_i)$ and $T' \in B(k' \varpi_{i'})$, if it exists, $T * T'$ denotes the tableau in $B(k \varpi_i + k' \varpi_{i'})$
obtained from $T$ by adding $T'$ to the right of $T$. Fix $b \in B(k \theta) \subset \Bad{\ell}$. Let
\begin{align*}
 b_L := \mathfrak{L}_k(b), \  b_R :=  \mathfrak{R}_k(b),\
 \varepsilon^1 := \nu_{1}(b_R), \  \varphi^1 := \nu_{n+1}(b_R), \
  \varepsilon^2 := \overline{\nu}_{n+1}(b_L), \  \varphi^2 := \overline{\nu}_{1}(b_L).
\end{align*}
The Kashiwara operators $\e_0, \f_0$ are given as follows:
\begin{align*}
\e_0(b) &:= \left\{
             \begin{array}{ll}
               b_L *  \e_0(b_R)  & \hbox{ if } \varphi^1 \ge \varepsilon^2 \text{ and } \varepsilon^1 > 0  , \\
               (b_L * \sBox{\overline{1}}\ ) * (b_R * \ssBox{n+1}\ ) & \hbox{ if } \varphi^1 \ge \varepsilon^2, \varepsilon^1 = 0 \text{ and } k < \ell \\
               \mathfrak{L}_1(\e_0(b_L))* \mathfrak{R}_1(b_R) & \hbox{ if } \varphi^1 < \varepsilon^2 \text{ and } \varepsilon^1 > 0 ,\\
               \e_0 (b_L) * b_R   & \hbox{ if }  \varphi^1 < \varepsilon^2  \text{ and } \varepsilon^1 = 0, \\
               0 & \hbox{otherwise,}
             \end{array}
           \right.\\
\f_0(b) &:= \left\{
             \begin{array}{ll}
               \mathfrak{L}_1(b_L) * \mathfrak{R}_1( \f_0(b_R) ) & \hbox{ if } \varphi^1 > \varepsilon^2 \text{ and } \varphi^2 > 0, \\
               b_L * \f_0(b_R) & \hbox{ if } \varphi^1 > \varepsilon^2 \text{ and } \varphi^2 = 0, \\
               \f_0(b_L)*b_R & \hbox{ if } \varphi^1 \le \varepsilon^2 \text{ and } \varphi^2 > 0 ,\\
               (\ \ssBox{\overline{n+1}} * b_L) * (\ \sBox{1} * b_R)   & \hbox{ if }  \varphi^1 \le \varepsilon^2 , \varphi^2 = 0 \text{ and } k < \ell, \\
               0 & \hbox{otherwise.}
             \end{array}
           \right.
\end{align*}
Then the crystal $\Bad{\ell}$ with $\e_0, \f_0$ becomes a $U_q'(\g)$-crystal.
\begin{Lem}\cite[Section 3.3]{Kde09} \label{Lem:B1 Bn isom Bad}
Let $ \xi^{\ell} : \Bone{\ell} \otimes \Ben{\ell} \rightarrow  \Bad{\ell}$ be the map defined by
$$ \xi^{\ell}(b\otimes \overline{b}) := \left\{
                                      \begin{array}{ll}
                                        \mathfrak{L}_{\overline{\nu}_1(\overline{b})}(\overline{b})*\mathfrak{R}_{\overline{\nu}_1(\overline{b})}(b) & \hbox{ if } \nu_1(b) \ge \overline{\nu}_1(\overline{b}). \\
                                        \mathfrak{L}_{\nu_1(b)}(\overline{b})*\mathfrak{R}_{\nu_1(b)}(b) & \hbox{ otherwise,}
                                      \end{array}
                                    \right.
 $$
for $b\in \Bone{\ell}$ and $\overline{b} \in \Ben{\ell}$. Then $\xi^{\ell}$ is a crystal isomorphism.
\end{Lem}
Note that $\Bad{\ell}$ is perfect by \cite[Lemma 4.6.2]{KKMMNN91}.

\begin{Lem} \label{Lem:g-s path in Bad}
Let $\Lambda := \sum_{k\in I} a_k\Lambda_k $ be a dominant integral weight of level $\ell$. We set $\mathtt{b}$ to be
the tableau in $B( ( \ell - a_0) \theta )$ such that the number of $i$ in $\mathtt{b}$ and the number of $\overline{i}$ in $\mathtt{b}$ are equal to $a_{i-1}$
for $2\le i \le n+1$.
Then
$$ \mathbf{p}^{\rm ad} (\Lambda) := \cdots \otimes \mathtt{b}\otimes \mathtt{b}\otimes \mathtt{b} $$
is the ground-state path of $\Lambda$ in $\Bad{\ell}$.
\end{Lem}
\begin{proof}
We write $ \Lambda = \Lambda_{j_1} + \cdots + \Lambda_{j_\ell} $ for $ j_1 \le \ldots \le j_\ell$. Define a tableau
$$
\mathtt{b} := \begin{tabular}{|c|c|c|c|c|c|}
          \hline
          $\overline{j_\ell +1}$ & $\cdots$ & $\overline{j_{a_0+1} +1}$ & $j_{a_0+1} +1$ & $\cdots$ & $j_\ell +1$ \\
          \hline
        \end{tabular} \ .
$$
By construction, $\mathtt{b}$ is contained in $B( ( \ell - a_0) \theta )$ and such that the number of $i$ in $\mathtt{b}$ and the number of $\overline{i}$ in $\mathtt{b}$ are equal to $a_{i-1}$
for $1\le i \le n+1$. Therefore, the assertion follows from $ \varepsilon(\mathtt{b}) = \varphi(\mathtt{b}) = \Lambda$.
\end{proof}

We denote by $\mathcal{P}^{\rm ad}(\Lambda)$ the path realization arising from $\Bad{\ell}$.
From $\eqref{Eq:map wt to B1}$, $\eqref{Eq:map wt to Bn}$ and Lemma \ref{Lem:B1 Bn isom Bad}, we obtain the 1-1 correspondence between $\wt(\Bone{\ell}) \times \wt(\Ben{\ell})$ and $\Bad{\ell}$:
\begin{equation} \label{Eq:map wt to Bad}
\begin{aligned}
\psi^{{\rm ad}, \ell}: \ & \wt(\Bone{\ell}) \times \wt(\Ben{\ell}) \longrightarrow \ \Bad{\ell},\\
& \qquad \quad \  (\alpha, \beta) \quad \qquad  \longmapsto\  \xi^{\ell}( \psi^{1,\ell}(\alpha) \otimes \psi^{n,\ell}(\beta) ).
\end{aligned}
\end{equation}

\vskip 1em

\begin{Exa} \label{Ex:perfect crystals}

Let $\g$ be of type $A_2^{(1)}$ and $\Lambda := 2\Lambda_0+\Lambda_1$.
We set $b := \f_1^4 \f_2^5 \f_1^2 \f_0^4 \f_2 \f_1 u_\Lambda$ in $B(\Lambda)$, where $u_\Lambda$ is the maximal vector of $B(\Lambda)$.
Then the corresponding elements of $b$ in $\mathcal{P}^1(\Lambda)$, $\mathcal{P}^n(\Lambda)$ and $\mathcal{P}^{\rm ad}(\Lambda)$ are
given as follows:
\begin{align*}
\mathbf{p}^1  &= \cdots \otimes
\begin{tabular}{|c|c|c|} \hline 2 & 2 & 3 \\ \hline \end{tabular} \otimes
\begin{tabular}{|c|c|c|} \hline 2 & 3 & 3 \\ \hline \end{tabular} \otimes
\begin{tabular}{|c|c|c|} \hline 2 & 2 & 3 \\ \hline \end{tabular} \otimes
\begin{tabular}{|c|c|c|} \hline 3 & 3 & 3 \\ \hline \end{tabular} \otimes
\begin{tabular}{|c|c|c|} \hline 1 & 2 & 3 \\ \hline \end{tabular}\ , \\
\mathbf{p}^n  &= \cdots \otimes
\begin{tabular}{|c|c|c|} \hline $\overline{3}$ & $\overline{3}$ & $\overline{1}$ \\ \hline \end{tabular} \otimes
\begin{tabular}{|c|c|c|} \hline $\overline{2}$ & $\overline{2}$ & $\overline{2}$ \\ \hline \end{tabular} \otimes
\begin{tabular}{|c|c|c|} \hline $\overline{1}$ & $\overline{1}$ & $\overline{1}$ \\ \hline \end{tabular} \otimes
\begin{tabular}{|c|c|c|} \hline $\overline{3}$ & $\overline{3}$ & $\overline{2}$ \\ \hline \end{tabular} \otimes
\begin{tabular}{|c|c|c|} \hline $\overline{1}$ & $\overline{1}$ & $\overline{1}$ \\ \hline \end{tabular}\ \otimes
\begin{tabular}{|c|c|c|} \hline $\overline{2}$ & $\overline{1}$ & $\overline{1}$ \\ \hline \end{tabular}\ , \\
\mathbf{p}^{\rm ad}  &= \cdots \otimes
\begin{tabular}{|c|c|} \hline 1 & 2 \\ \cline{1-2} 3  \\ \cline{1-1} \end{tabular} \otimes
\begin{tabular}{|c|c|} \hline 1 & 2 \\ \cline{1-2} 3  \\ \cline{1-1} \end{tabular} \otimes
\begin{tabular}{|c|c|} \hline 2 & 3 \\ \cline{1-2} 3  \\ \cline{1-1} \end{tabular}\otimes
\begin{tabular}{|c|c|c|c|c|c|}
  \hline
  1 & 2 & 2 & 2 & 2 & 3 \\
 \hline
  3 & 3 & 3  \\
  \cline{1-3}
\end{tabular} \ .
\end{align*}

\end{Exa}

\vskip 2em

\section{Geometric realizations of crystals}

In this section, we review the geometric constructions $\gB(\Lambda)$ for type $A_n^{(1)}$.
Let $I = \Z / (n+1)\Z$ and $H$ the set of the arrows such that $i \rightarrow j$ with
$i,j \in I,\ i-j= \pm 1\ $. For $h \in H$, we denote by
$\qin(h)$ (resp.\ $\qout(h)$) the incoming (resp.\ outgoing) vertex
of $h$. Define an involution $^-:H \to H$ to be the map
interchanging $i\to j$ and $j\to i$. Let
$\Omega := \{ h \in H \mid \qin(h) - \qout(h)  =1   \}$
so that $H = \Omega \sqcup
\overline{\Omega}$; i.e.,

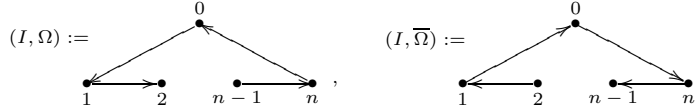
$$ \Quiver$$

We take the map $\epsilon:H \to \{-1,1\}$ given by $\epsilon(h) := 1$ if $h \in \Omega$ and $\epsilon(h) := -1$ if $h \in \overline{\Omega}$.
Given an $I$-graded vector space $V = \bigoplus_{i=0}^n V_i$, let $\underline{\dim} V := \sum_{i=0}^n \dim(V_i) \alpha_i$.
For $\alpha = \sum_{i=0}^n k_i \alpha_i \in Q^+$, let $V(\alpha) := \bigoplus_{i=0}^n V_i(\alpha)$ be an $I$-graded vector space
such that $V_i(\alpha)$ is a $\C$-vector space with an ordered basis $v^i(\alpha) = \{ v^i_0,v^i_1 \ldots, v^i_{k_i-1} \}$ for $i\in I$.
Fix an ordered basis
$$v(\alpha)=\{ v^0_0, \ldots, v^0_{k_0-1}, v^1_{0}, \ldots, v^1_{k_1-1}, \ldots, v^n_0, \ldots  v^n_{k_n-1} \},$$
for $V(\alpha)$.
For simplicity, we write $V$ (resp.\ $V_i$) for $V(\alpha)$ (resp.\ $V_i(\alpha)$). Let
$$ E(\alpha) := E_\Omega (\alpha ) \oplus E_{\overline{\Omega}}(\alpha),$$
where $E_{\Omega}(\alpha ) := \bigoplus_{i \in I}\Hom(V_{i-1},\ V_{i}) $ and $ E_{\overline{\Omega}}(\alpha) :=  \bigoplus_{i \in I}\Hom(V_{i},\ V_{i-1}).$

Let $\pi_{\Omega}$ (resp.\ $\pi_{\overline{\Omega}}$) denote the natural projection from $E(\alpha)$ to $E_{\Omega}(\alpha)$ (resp.\ $E_{\overline{\Omega}}(\alpha)$).
For an element $\chi \in E(\alpha)$, if there is no danger of confusion, we write $x=(x_i\in \Hom(V_{i-1},\ V_{i}))_{i\in I}$ (resp.\ $\overline x=( \overline{x}_i \in \Hom(V_{i},\ V_{i-1}))_{i\in I}$)
 for $\pi_{\Omega}(\chi)$ (resp.\ $\pi_{\overline \Omega}(\chi)$). The matrix representations of $x$ and $\overline{x}$ in the ordered basis $v(\alpha)$ are given as
$$ x = \left(
         \begin{array}{ccccc}
           0   &  &   \cdots  &       0  & x_0 \\
           x_1  & 0   &  &         &  0 \\
           0   &   x_2  & 0   &   & \vdots \\
           \vdots    &     &  \ddots   & \ddots  & 0 \\
           0 &  \cdots   &  0   &     x_n    & 0 \\
         \end{array}
       \right),
       \qquad
       \overline{x} = \left(
         \begin{array}{ccccc}
           0   & \overline{x}_1 & 0    &     \cdots  & 0 \\
            & 0   & \overline{x}_2 &         &   \vdots \\
            \vdots   &     & 0   &  \ddots & 0 \\
             0  &     &    & \ddots  & \overline{x}_n  \\
           \overline{x}_0 &  0   &  \cdots   &     0    & 0 \\
         \end{array}
       \right),
 $$
where $x_i$ (resp.\ $\overline{x}_i$) is the matrix representation of $x|_{V_{i-1}}: V_{i-1} \to V_i$ (resp.\ $\overline{x}|_{V_i}: V_i \to V_{i-1}$)
in the ordered bases $v^{i-1}(\alpha)$ and $v^{i}(\alpha)$ (resp.\ $v^i(\alpha)$ and $v^{i-1}(\alpha)$).

The algebraic group $G(\alpha):= \prod_{i \in I} \Aut(V_i) \subset \Aut(V)$ acts on $E(\alpha)$ by
$g \cdot\chi  := g \chi g^{-1}$ for $g \in G(\alpha)$ and $ \chi \in E(\alpha)$.
For $ \chi \in E(\alpha)$ and $h\in H$, define $\chi_{h} := x_{\qin(h)}$ if $h\in \Omega$ and $\chi_{h} := x_{\qout(h)}$ if $h\in \overline{\Omega}$.
Let $\langle \ ,\ \rangle$ be the sympletic form on $E(\alpha)$ defined by
$ \langle \chi,\chi' \rangle := \sum_{h\in H} \epsilon(h) \tr(\chi_h {\chi'}_{\overline h})$
for $\chi,\chi' \in E(\alpha)$. Note that $E(\alpha)$ may be viewed as the cotangent bundle of $E_{\Omega}(\alpha)$ (resp.\ $E_{\overline{\Omega}}(\alpha)$) under this form. The {\it moment map} $\mu=(\mu_i: E(\alpha) \to \End(V_i))_{i\in I}$ is given by
$ \mu_i(\chi) := \sum_{h\in H,\ \qin(h)=i} \epsilon(h)\chi_h \chi_{\overline h}$
for $\chi \in E(\alpha)$. Note that, for $\chi = x + \overline{x}\in E(\alpha)$,
\begin{align}\label{Eq:commuting property of chi}
\mu_i(\chi)=0\ \text{for all } i \in I\quad \text{ if and only if } \quad x \overline x - \overline x  x = 0.
\end{align}

An element $\chi \in E(\alpha)$ is {\it nilpotent} if there exists
an $N \ge 2$ such that for any sequence $h_1, \ldots, h_N \in H$
satisfying $\qin(h_i)=\qout(h_{i+1})\ (i=1,\ldots,N-1)$, the
composition map $\chi_{h_N} \cdots \chi_{h_1}$ is zero. We define
{\it Lusztig's quiver variety} \cite{L90, L91} to be
$$\mathfrak{X}(\alpha) := \{\chi \in E(\alpha) \mid \chi: \text{nilpotent},\  \mu_i(\chi)=0 \text{ for all }i\in I \}.$$
Let $\Irr \mathfrak{X}(\alpha)$ denote the set of all irreducible components of $\mathfrak{X}(\alpha)$.
In \cite{KS97}, Kashiwara and Saito gave a crystal structure to
$ \gB(\infty) := \bigsqcup_{\alpha \in Q^+} \Irr \mathfrak{X} (\alpha),$
and proved the following theorem.
\begin{Thm}\cite[Theorem 5.3.2]{KS97}
There is a unique crystal isomorphism $\gB(\infty) \simeq B(\infty)$.
\end{Thm}

We now review Nakajima's quiver varieties \cite{N94,N98}. For $\Lambda = \sum_{i=0}^n a_i \Lambda_i \in P^+ $, let
$ W(\Lambda) := \bigoplus_{i=0}^n W_i(\Lambda)$ such that $W_i(\Lambda)$ is a $\C$-vector space of dimension $a_i$.
For simplicity, we write $W$ (resp.\ $W_i$) for $W(\Lambda)$ (resp.\ $W_i(\Lambda)$).
Let
$$ X(\Lambda, \alpha) := \mathfrak{X}(\alpha) \times \sum_{i \in I} \Hom(V_i, W_i). $$
The group $G(\alpha)$ acts on $X(\Lambda, \alpha)$ by $g \cdot(\chi,
t) := (g \chi g^{-1}, t g^{-1})$. For $\chi
\in \mathfrak{X}(\alpha)$, an $I$-graded subspace $S$ of $V(\alpha)$ is
{\it $\chi$-invariant} if $ \chi_{h}(S_{\qout(h)}) \subset
S_{\qin(h)} $ for all $h \in H $. An element $(\chi, t)
\in X(\Lambda, \alpha )$ is called a {\it stable point} if it satisfies the following conditions: if
$S$ is a $\chi$-invariant subspace of $V$ with $t_i(S_i)=0\
(i\in I)$, then $S = 0$. For $(\chi = x+ \overline{x}, t)\in  X(\Lambda, \alpha )$,
one can prove that
\begin{align} \label{Eq:equiv of stable pts}
(\chi, t) \text{ is stable }  \Longleftrightarrow\ \ker x \cap \ker \overline{x} \cap \ker t = 0.
\end{align}
Let $X(\Lambda, \alpha)^{\rm st}$ be the set
of all stable points of $X(\Lambda, \alpha)$. We define
$$\mathfrak{X}(\Lambda, \alpha) := X(\Lambda, \alpha)^{\rm st} / G(\alpha).$$
Let $\Irr \mathfrak{X}(\Lambda, \alpha)$ (resp.\ $\Irr X(\Lambda, \alpha)$) be the set of all irreducible components
of $\mathfrak{X}(\Lambda, \alpha)$ (resp.\ $X(\Lambda, \alpha)$).
Since $\Irr \mathfrak{X}(\Lambda, \alpha)$ can be identified with
$ \{ Y \in \Irr X(\Lambda, \alpha)|\ Y \cap X(\Lambda, \alpha)^{\rm st} \ne \emptyset \}, $
we have
\begin{equation} \label{Eq:def of Pi}
\begin{aligned}
\Pi_\Lambda :  & \  \Irr \mathfrak{X}(\alpha) \  \longrightarrow \  \Irr \mathfrak{X}(\Lambda, \alpha) \sqcup \{ 0 \}, \\
& \quad \   X_0 \quad \longmapsto \quad \left( \left( X_0 \times \sum_{i\in I} \Hom(V_i, W_i)  \right) \cap X(\Lambda, \alpha)^{\rm st} \right) / G(\alpha).
\end{aligned}
\end{equation}
In \cite{St02}, Saito defined a crystal structure on
$ \gB(\Lambda) := \bigsqcup_{\alpha \in Q^+} \Irr \mathfrak{X}(\Lambda, \alpha) ,$
and proved the following theorem.

\begin{Thm}\cite[Theorem 4.6.4]{St02} \label{Thm:saito geometric crystal}
There is a unique crystal isomorphism $\gB(\Lambda) \simeq B(\Lambda)$.
\end{Thm}

\vskip 2em

\section{Quiver varieties and Young walls}

The {\em Young walls} are new combinatorial objects introduced by Kang \cite{HK02,Kang03} for realizing crystals of basic representations for several types.
Young walls generalize the notion of colored Young diagrams in \cite{DJKMO89, MM90}.
In this section, we deduce two kinds of Young walls from the results given in \cite{FLOTW99, JMMO91, KL06, Sv06},
and describe irreducible components of $\gB(\Lambda)$ by using these two kinds of Young walls.
This description will play a crucial role in proving Theorem \ref{Thm:Fundamantal thm for qv}.

Let $\mathtt{P}_k^1$ and $\mathtt{P}_k^n$ be the patterns of a fundamental weight $\Lambda_k$ given as follows:

$$ \Patterns $$

A Young wall $Y$ in $\mathtt{P}_k^1$ (resp.\ $\mathtt{P}_k^n$) is a wall consisting of colored blocks stacked by the following rules:
\begin{enumerate}
\item[(a)] the colored blocks should be stacked in the pattern $\mathtt{P}_k^1$ (resp.\ $\mathtt{P}_k^n$) of weight $\Lambda_k$,
\item[(b)] except for the rightmost column, there should be no free space to the right of any block.
\end{enumerate}
We number columns of $Y$ from right to left,
denoted by $Y = (y_j)_{j \ge 0}$ where $y_j$ is the $j$th column of $Y$. Let $\HT(y_j)$ be the number of blocks in $y_j$, and
$\wt(y_j) := \sum_{i\in I} k_{ij} \alpha_i\ (j \in \Z_{\ge0} ),$
where $k_{ij}$ is the number of $i$-blocks in the $j$th column $y_j$. Let $\mathcal{Y}^1_k$ (resp.\ $\mathcal{Y}^n_k$) denote the set of all Young walls
in the pattern $\mathtt{P}_k^1$ (resp.\ $\mathtt{P}_k^n$).
Fix $\Lambda \in P^+$ of level $\ell$. Then $\Lambda$ can be uniquely written as
$$ \Lambda = \Lambda_{i_1} + \cdots + \Lambda_{i_\ell}$$
for some $0 \le i_1 \le i_2 \le \cdots \le i_\ell \le n$.
For an $\ell$-tuple $\mathbf{Y} = (Y^{(1)}, \ldots, Y^{(\ell)} )$ of Young walls, let $l^{(k)}_j$ be the length of the $j$th row of $Y^{(k)}$,
and $\col^{(k)}_j$ be the color at the left end of the $j$th row of $Y^{(k)}$.
An $\ell$-tuple $\mathbf{Y} = (Y^{(1)}, \ldots, Y^{(\ell)} ) $ of Young walls is said to be {\em reduced} if, for all $t > 0$,
$$ \{ \col^{(k)}_j \mid l^{(k)}_j = t,\ 1 \le k \le \ell,\ j \ge 0 \} \ \neq \   \{0,1,\ldots, n \} .$$
Let
\begin{equation} \label{Eq:def of Y1}
\begin{aligned}
\mathcal{Y}^1(\Lambda) := \{ \mathbf{Y} = (Y^{(1)}, \ldots, Y^{(\ell)} ) \mid  &\  Y^{(k)} \in \mathcal{Y}^1_{i_k},\ \mathbf{Y}: \text{reduced}, \\
 & Y^{(1)} \preceq Y^{(2)} \preceq \cdots \preceq Y^{(\ell)}\preceq Y^{(1)}+n+1 \} .
\end{aligned}
\end{equation}
Here, if let $Y^{(k)} = (y^{(k)}_j)_{j \ge 0}$, then the order $\preceq$ is given as follows:
\begin{align*}
Y^{(k)} \preceq Y^{(k+1)} \ &\Longleftrightarrow \ \HT(y^{(k)}_j) \le \HT(y^{(k+1)}_j) + i_{k+1}-i_k \ (j \in \Z_{\ge 0}), \\
Y^{(\ell)} \preceq Y^{(1)} +n+1   \ &\Longleftrightarrow \ \HT(y^{(\ell)}_j) \le \HT(y^{(1)}_j) +i_1 - i_{\ell} + n+1 \ (j \in \Z_{\ge 0}).
\end{align*}
Then, by \cite[Proposition 2.11]{FLOTW99}, we have a 1-1 map
\begin{equation} \label{Eq:iso between Y1 and P1}
\begin{aligned}
& F^1 \  : \ \qquad  \mathcal{Y}^1(\Lambda) \qquad \ \   \longrightarrow \quad \mathcal{P}^1(\Lambda), \\
& \ \  \mathbf{Y} = (Y^{(1)}, \ldots, Y^{(\ell)} ) \quad \longmapsto \   \left( \psi^{1,\ell} \left(\wt(b_j^\Lambda) -\cl \left( \sum_{k=1}^\ell \wt(y_j^{(k)}) \right) \right) \right)_{j \ge 0 },
\end{aligned}
\end{equation}
where $\mathbf{p}^1(\Lambda) = ( b^{\Lambda}_k )_{k \ge 0}$ is the ground-state path of $\Lambda$ in $\Bone{\ell}$.

In a similar manner, we let
\begin{equation}  \label{Eq:def of Yn}
\begin{aligned}
\mathcal{Y}^n(\Lambda) := \{ \overline{\mathbf{Y}} = (\overline{Y}^{(1)}, \ldots, \overline{Y}^{(\ell)} ) \mid &\ \overline{Y}^{(k)} \in \mathcal{Y}^n_{i_k},\ \overline{\mathbf{Y}}: \text{reduced}, \\
&\ \overline{Y}^{(1)} \succeq\overline{Y}^{(2)} \succeq \cdots \succeq \overline{Y}^{(\ell)}\succeq \overline{Y}^{(1)}-n-1 \} .
\end{aligned}
\end{equation}
Here, if let $\overline{Y}^{(k)} = (\overline{y}^{(k)}_j)_{j \ge 0}$, then  the order $\succeq$ is given as follows:
\begin{align*}
\overline{Y}^{(k-1)} \succeq \overline{Y}^{(k)} \ &\Longleftrightarrow \ \HT(\overline{y}^{(k-1)}_j)  \ge \HT(\overline{y}^{(k)}_j) + i_{k-1} - i_k   \ (j \in \Z_{\ge 0}), \\
\overline{Y}^{(\ell)} \succeq \overline{Y}^{(1)} -n-1   \ &\Longleftrightarrow \ \HT(\overline{y}^{(\ell)}_j) \ge \HT(\overline{y}^{(1)}_j) + i_{\ell} - i_1- n-1 \ (j \in \Z_{\ge 0}).
\end{align*}
We have a 1-1 map
\begin{equation} \label{Eq:iso between Yn and Pn}
\begin{aligned}
& F^n \  : \ \qquad  \mathcal{Y}^n(\Lambda) \qquad \ \   \longrightarrow \quad \mathcal{P}^n(\Lambda), \\
& \ \  \overline{\mathbf{Y}} = (\overline{Y}^{(1)}, \ldots, \overline{Y}^{(\ell)} ) \quad \longmapsto \   \left(  \psi^{n,\ell} \left(\wt(\overline{b}_j^\Lambda) -\cl \left( \sum_{k=1}^\ell \wt(\overline{y}_j^{(k)}) \right) \right)  \right)_{j \ge 0 },
\end{aligned}
\end{equation}
where $\mathbf{p}^n(\Lambda) = ( \overline{b}^{\Lambda}_k )_{k \ge 0}$ is the ground-state path of $\Lambda$ in $\Ben{\ell}$.
Then, $\mathcal{Y}^1(\Lambda)$ and $\mathcal{Y}^n(\Lambda)$ have the crystal structures from $\mathcal{P}^1(\Lambda)$ and $\mathcal{P}^n(\Lambda)$ via the maps $F^1$ and $F^n$, respectively.

\begin{Prop} \label{Prop:iso of YW and Path}
 For $\Lambda \in P^+$, $\mathcal{Y}^1(\Lambda)$ (resp.\ $\mathcal{Y}^n(\Lambda)$) is isomorphic to
the crystal $B(\Lambda)$.
\end{Prop}
Note that the crystal structures of $\mathcal{Y}^1(\Lambda)$ and $\mathcal{Y}^n(\Lambda)$ can be described completely in a combinatorial manner
by using \cite[Proposition 2.11]{FLOTW99}, \cite[Theorem 3.13]{JMMO91} and \cite[Proposition 6.1]{KL06}.
As a set, $\mathcal{Y}^1(\Lambda)$ is different from $\mathcal{Y}^n(\Lambda)$ even when forgetting colors of Young walls. It is difficult to describe
an explicit crystal isomorphism between $\mathcal{Y}^1(\Lambda)$ and $\mathcal{Y}^n(\Lambda)$ in a combinatorial manner \cite{JL09}.
We also remark that the description of $\mathcal{Y}^n(\Lambda)$ is almost the same as that of {\it Young pyramid} given in \cite{Sv06}.
Given $\overline{\mathbf{Y}} = (\overline{Y}^{(1)}, \ldots, \overline{Y}^{(\ell)} ) \in \mathcal{Y}^n(\Lambda)$,
one can obtain the corresponding Young pyramid by flipping $\overline{Y}^{(k)}$ in the diagonal and rearranging them in an appropriate manner.

Now we describe irreducible components of quiver varieties using $\mathcal{Y}^1(\Lambda)$ and $\mathcal{Y}^n(\Lambda)$.
Take $\alpha := \sum_{i=0}^n k_i \alpha_i \in Q^+$ and write $V$ (resp.\ $V_i$) for $V(\alpha)$ (resp.\ $V_i(\alpha)$). Let $k_{-1}=k_n$ and $k_{n+1}=k_{0}$.
By Proposition \ref{Prop:iso of YW and Path} and Theorem \ref{Thm:saito geometric crystal}, there is the unique crystal isomorphism from $\mathcal{Y}^1(\Lambda)$
(resp.\ $\mathcal{Y}^n(\Lambda)$) to $\gB(\Lambda)$. We will describe these crystal isomorphisms explicitly.

For $s \in I$, we define a linear map $\mathcal{E}^s_{ij}:V_{s-1} \to V_{s}$ (resp.\ $\overline{\mathcal{E}}^s_{ij}:V_{s} \to V_{s-1}$) by
$$
 \mathcal{E}^{s}_{ij}(v^{s-1}_k) := \left\{
                       \begin{array}{ll}
                         v^{s}_j & \hbox{ if } k=i, \\
                         0 & \hbox{ otherwise,}
                       \end{array}
                   \right.
                   \qquad \left( \text{resp.\ } \ \
\overline{\mathcal{E}}^{s}_{ij}(v^{s}_k) := \left\{
                       \begin{array}{ll}
                         v^{s-1}_j & \hbox{ if } k=i, \\
                         0 & \hbox{ otherwise.}
                       \end{array}
                   \right.
                   \right)
$$
We first deal with an isomorphism of $\mathcal{Y}^n(\Lambda)$ and $ \gB(\Lambda) $.
Take $\overline{\mathbf{Y}} = (\overline{Y}^{(1)}, \ldots, \overline{Y}^{(\ell)}) \in \mathcal{Y}^n(\Lambda)$. Let $\overline{\mathfrak{b}}_{ij}^{(k)}$ denote the $i$th
block from the bottom of the $j$th column of $\overline{Y}^{(k)}$ and $\col(\overline{\mathfrak{b}}_{ij}^{(k)})$ be the color of $\overline{\mathfrak{b}}_{ij}^{(k)}$.
We set
$ o(\overline{\mathfrak{b}}_{ij}^{(k)}) := \# \{ \overline{\mathfrak{b}}_{rs}^{(t)} \in \overline{\mathbf{Y}} \mid
 \col(\overline{\mathfrak{b}}_{ij}^{(k)}) = \col(\overline{\mathfrak{b}}_{rs}^{(t)}),\ (t,r,s) \prec (k,i,j) \}, $
where $\prec$ is the lexicographical order. Define
\begin{align} \label{Eq:def of bar x(bar Y)}
\overline{x}(\overline{\mathbf{Y}}) := \sum_{\overline{\mathfrak{b}}_{ij}^{(k)} \in \overline{\mathbf{Y}}, j > 0}
\overline{\mathcal{E}}^{\col(\overline{\mathfrak{b}}_{ij}^{(k)})}_{o(\overline{\mathfrak{b}}_{ij}^{(k)}),\ o(\overline{\mathfrak{b}}_{ij-1}^{(k)})} \in E_{\overline{\Omega}}(\alpha).
\end{align}
Note that, for $t \ge 0$,
\begin{align} \label{Eq:computation of ker bar Y}
\ker (\overline{x}(\overline{\mathbf{Y}}))^t = \Span_\C\{ v_{o(\overline{\mathfrak{b}}_{ij}^{(k)}) } ^{\col(\overline{\mathfrak{b}}_{ij}^{(k)})}
 \mid \overline{\mathfrak{b}}_{ij}^{(k)} \in \overline{Y}^{(k)},\ j < t,\ 1\le k \le \ell  \}.
\end{align}
Since $E(\alpha)$ can be viewed as the cotangent bundle of $E_{\overline{\Omega}}(\alpha)$, we let
\begin{align} \label{Eq:def of irr component of bar Y}
X_0(\overline{\mathbf{Y}}) := \text{the closure of the conormal bundle of the $G(\alpha)$-orbit of $\overline{x}(\overline{\mathbf{Y}})$}.
\end{align}
\begin{Prop} \label{Prop:isom between Yn and gB}
The map $\overline{\mathbf{Y}} \mapsto \Pi_\Lambda(X_0(\overline{\mathbf{Y}}))$ from $\mathcal{Y}^n(\Lambda)$ to $ \gB(\Lambda) $ is a crystal isomorphism.
\end{Prop}
\begin{proof}
By the definition of Young walls, the pair $(\overline{x}(\overline{\mathbf{Y}}), V(\alpha))$ becomes a nilpotent representation of the quiver $(I, \overline{\Omega})$. Moreover,
it follows from the definition of $\mathcal{Y}^n(\Lambda)$ in $\eqref{Eq:def of Yn}$ that $(\overline{x}(\overline{\mathbf{Y}}), V(\alpha))$
satisfies the {\em aperiodic condition} given in \cite[Section 15.4]{L91}.
Hence, by \cite[Corollary 15.6]{L91}, $X_0(\overline{\mathbf{Y}})$ is an irreducible component of $\mathfrak{X}(\alpha)$.
Therefore, since the description of $\mathcal{Y}^n(\Lambda)$ is almost the same as that of Young pyramid given in \cite{Sv06},
our assertion follows from \cite[Theorem 6.3]{FS03} and \cite[Theorem 8.4]{Sv06}.
\end{proof}

Now we consider the crystal isomorphism between $\mathcal{Y}^1(\Lambda)$ and $ \gB(\Lambda) $.
Let $\mathbf{Y} = (Y^{(1)}, \ldots, Y^{(\ell)}) \in \mathcal{Y}^1(\Lambda)$. Let $\mathfrak{b}_{ij}^{(k)}$ denote the $i$th
block from the bottom of the $j$th column of $Y^{(k)}$, and $\col(\mathfrak{b}_{ij}^{(k)})$ be the color of $\mathfrak{b}_{ij}^{(k)}$. We set
$ o(\mathfrak{b}_{ij}^{(k)}) := \# \{ \mathfrak{b}_{rs}^{(t)} \in \mathbf{Y} \mid
 \col(\mathfrak{b}_{ij}^{(k)}) = \col(\mathfrak{b}_{rs}^{(t)}),\ (t,r,s) \prec (k,i,j) \}, $
where $\prec$ is the lexicographical order. Define
\begin{align} \label{Eq:def of x(Y)}
x(\mathbf{Y}) := \sum_{\mathfrak{b}_{ij}^{(k)} \in \mathbf{Y}, j > 0}
\mathcal{E}^{\col(\mathfrak{b}_{i,j-1}^{(k)})}_{o(\mathfrak{b}_{ij}^{(k)}),\ o(\mathfrak{b}_{ij-1}^{(k)})} \in E_\Omega(\alpha).
\end{align}
Note that, for $t \ge 0$,
\begin{align} \label{Eq:computation of ker Y}
\ker (x(\mathbf{Y}))^t = \Span_\C\{ v_{o(\mathfrak{b}_{ij}^{(k)}) } ^{\col(\mathfrak{b}_{ij}^{(k)})} \mid \mathfrak{b}_{ij}^{(k)} \in Y^{(k)},\ j < t,\ 1\le k \le \ell  \}.
\end{align}
Since $E(\alpha)$ can be viewed as the cotangent bundle of $E_{\Omega}(\alpha)$, let
\begin{align} \label{Eq:def of irr component of Y}
X_0(\mathbf{Y}) := \text{the closure of the conormal bundle of the $G(\alpha)$-orbit of $x(\mathbf{Y})$}.
\end{align}

\begin{Prop} \label{Prop:isom between Y1 and gB}
The map $\mathbf{Y} \mapsto \Pi_\Lambda(X_0(\mathbf{Y}))$ from $\mathcal{Y}^1(\Lambda)$ to $ \gB(\Lambda) $ is a crystal isomorphism.
\end{Prop}
\begin{proof}
Note that $\mathcal{Y}^1(\Lambda)$ and $\mathcal{Y}^n(\Lambda)$ are totally different as sets.
By the definition of Young walls, the pair $(x(\mathbf{Y}), V(\alpha))$ becomes a nilpotent representation of the quiver $(I, \Omega)$.
Since the counterclockwise orientation $\Omega$ is opposite to $\overline{\Omega}$,
by the symmetry of $\Omega$ and $\overline{\Omega}$, one can prove that the pair $(x(\mathbf{Y}), V(\alpha))$ satisfies the aperiodic condition
and $X_0(\mathbf{Y})$ is an irreducible component of $\mathfrak{X}(\alpha)$.
For the same reason, our assertion can be proved in a similar manner to that of \cite[Theorem 8.4]{Sv06}.
\end{proof}

\vskip 1em

\begin{Exa} \label{Ex:Young walls}
We keep all notations in Example \ref{Ex:perfect crystals}. The Young walls corresponding to $b$ in $\mathcal{Y}^1(\Lambda)$ and $\mathcal{Y}^n(\Lambda)$
are given as follows.
\begin{align*}
 \mathbf{Y} &= \left( \ \ \ExYWone \ \ \right), \\
 \overline{\mathbf{Y}} &= \left( \ \ \ExYWen \ \ \right).
\end{align*}

We set $ \alpha := \Lambda - \wt(b) = 4 \alpha_0 + 7 \alpha_1 + 6 \alpha_2. $
Let $X$ be the irreducible component of $\gB(\Lambda)$ corresponding to $b$ via the isomorphism in Theorem \ref{Thm:saito geometric crystal}.
The matrices $x(\mathbf{Y})$ and $x(\overline{\mathbf{Y}})$ are given as follows.
\begin{align*}
x(\mathbf{Y}) &:= \mathcal{E}^{0}_{00} + \mathcal{E}^{0}_{11} + \mathcal{E}^{0}_{32} + \mathcal{E}^{0}_{53}
+ \mathcal{E}^{1}_{24} + \mathcal{E}^{2}_{00} + \mathcal{E}^{2}_{21} + \mathcal{E}^{2}_{53} + \mathcal{E}^{2}_{64}, \\
\overline{x}(\overline{\mathbf{Y}}) &:= \overline{\mathcal{E}}^{0}_{23} + \overline{\mathcal{E}}^{1}_{00} + \overline{\mathcal{E}}^{1}_{21}
+ \overline{\mathcal{E}}^{1}_{52} + \overline{\mathcal{E}}^{1}_{63} + \overline{\mathcal{E}}^{2}_{00}
+ \overline{\mathcal{E}}^{2}_{34} + \overline{\mathcal{E}}^{2}_{45}.
\end{align*}
Note that $ x(\mathbf{Y}) + \overline{x}(\overline{\mathbf{Y}}) \notin \mathfrak{X}(\alpha) $
since $[x(\mathbf{Y}), \overline{x}(\overline{\mathbf{Y}})] \ne 0$.

\end{Exa}

\vskip 2em

\section{Quiver varieties and the crystals $B^{1,\ell}$, $B^{n,\ell}$ and $\Bad{\ell}$ }

In this section, we construct an explicit crystal isomorphism from the geometric realization $\gB(\Lambda)$ to the path realization
$\mathcal{P}^{\rm ad}(\Lambda)$ (resp.\ $\mathcal{P}^1(\Lambda)$, $\mathcal{P}^n(\Lambda)$) arising from $\Bad{\ell}$ (resp.\ $\Bone{\ell}$, $\Ben{\ell}$).
Let $X \in \gB(\Lambda) $ for $\Lambda \in P^+$ of level $\ell$
and $\alpha := \Lambda - \wt(X)$.
For an element $(\chi, t) \in X(\Lambda, \alpha)$, let $[\chi, t]$ denote the $G(\alpha)$-orbit of $(\chi, t)$ in $X(\Lambda, \alpha)$.
For an element $[\chi=x+\overline{x}, t] $ in a certain open subset of $X$, we will give an explicit description of the $\Lambda$-path in $\Bad{\ell}$ (resp.\ $\Bone{\ell}$, $\Ben{\ell}$) corresponding to
$X$ in terms of dimension vectors of $\ker (x\overline{x})^{k+1}/ \ker \overline{x}(x\overline{x})^k$ and $\ker \overline{x}(x\overline{x})^{k}/ \ker (x\overline{x})^k$
(resp.\ $\ker x^{k+1}/ \ker x^k$, $\ker \overline{x}^{k+1}/ \ker \overline{x}^k$). For this purpose, we need several lemmas.
Combining $\eqref{Eq:def of Pi}$ with results in \cite{KP10}, we have the following lemmas.

\begin{Lem} \cite[Lemma 5.1]{KP10} \label{Lem:chi statble}
Let $X \in \gB(\Lambda)$. For any element $[\chi= x+ \overline{x}, t] \in X $ and $k \in \Z_{\ge 0}$, we have
\begin{enumerate}
\item $\ker(x \overline{x})^k = \ker(\overline{x} x)^k $,
\item $\ker x^k $ and $ \ker\overline{x}^k $ are $\chi$-invariant,
\item $\ker(x \overline{x})^k$ is $\chi$-invariant.
\end{enumerate}
\end{Lem}

\begin{Lem}  \cite[Lemma 5.2]{KP10} \label{Lem:open set U in B1 and Bn}
Let $X \in \gB(\Lambda)$. There is a nonempty open subset $U \subset X$ such that
$$ \ker x^k \simeq \ker x^{\prime k}, \quad  \ker \overline{x}^k \simeq \ker \overline{x}^{\prime k} $$
for any elements $[\chi = x + \overline{x}, t]$, $[\chi' = x' + \overline{x}', t'] \in U$ and $k \in \Z_{\ge 0}$.
\end{Lem}


Let $M := \bigoplus_{i\in I } M_i$ be an $I$-graded vector space and $x \in \Hom(M, M)$ be a homogeneous linear map of degree $k$.
For an $x$-invariant $I$-graded subspace $N$ of $M$, let $x|_{M/N}$ (resp.\ $x|_N$) denote the homogeneous linear map of degree $k$ in $\Hom(M/N, M/N)$
(resp.\ $ \Hom(N, N)$) obtained by restricting $x$ to $M/N$ (resp.\ $N$).

\begin{Lem} \label{Lem:existance of the open set}
Let $M = \bigoplus_{i\in I} M_i$ and $L = \bigoplus_{i\in I} L_i$.
Let $x \in \bigoplus_{i\in I} \Hom(M_{i-1}, M_i)$, and let
$$ N:= \ker x \ \text{ and } \ x' := x |_{M / N}. $$
Take elements
$$\overline{x}' \in \bigoplus_{i\in I} \Hom(M_i / N_i, M_{i-1}/N_{i-1}) \  \text{ and }\  t' \in \bigoplus_{i\in I} \Hom( M_i/ N_i, L_{i+1})$$
such that
$$ [x', \overline{x}']=0\ \text{ and }\ \ker x' \cap \ker \overline{x}' \cap \ker t' = 0, $$
where $N_i$ is the $i$-subspace of $N$ for $i\in I$. Then there exist
$$\overline{x} \in \bigoplus_{i\in I} \Hom(M_i , M_{i-1}) \  \text{ and }\  t \in \bigoplus_{i\in I} \Hom( M_i, L_i)$$
such that
$$ [x,\overline{x}]=0, \ \ \overline{x}|_{M/N} = \overline{x}', \ \ t' = t\circ x |_{M/N}  \  \text{ and }\  \ker x \cap \ker \overline{x} \cap \ker t = 0. $$
\end{Lem}

\begin{proof}
By definition, $x|_{\ker x^2 / \ker x}$ gives an $I$-graded isomorphism between $ \ker x^2 / \ker x$ and $ \ker x \cap \im\ x$. We set
\begin{align*}
\ell &:= \dim L, \\
r &:= \dim(\ker x \cap \im\ x) = \dim (\ker x^2 / \ker x), \\
s &:= \dim( \ker x / \ker x \cap \im\ x ), \\
t &:= \dim(M / \ker x^2).
\end{align*}
Take a homogeneous ordered basis
$$ \omega := \{ \omega^{1}_1,\ldots, \omega^{1}_r, \omega^{2}_1,\ldots, \omega^{2}_s, \omega^{3}_1,\ldots, \omega^{3}_r, \omega^{4}_1,\ldots, \omega^{4}_t  \} $$
of $M$ such that
\begin{align*}
\ker x \cap \im\ x &= \Span_\C\{ \omega^{1}_1,\ldots, \omega^{1}_r  \}, \\
\ker x &= \Span_\C\{ \omega^{1}_1,\ldots, \omega^{1}_r , \omega^{2}_1,\ldots, \omega^{2}_s \}, \\
\ker x^2  &= \Span_\C\{ \omega^{1}_1,\ldots, \omega^{1}_r , \omega^{2}_1,\ldots, \omega^{2}_s, \omega^{3}_1,\ldots, \omega^{3}_r \}.
\end{align*}
Then the matrix representative of $x$ with respect to the basis $\omega$ is given as follows:
$$ x = \left(
         \begin{array}{cccc}
           0 & 0 & a & b \\
           0 & 0 & 0 & c \\
           0 & 0 & 0 & d \\
           0 & 0 & 0 & e \\
         \end{array}
       \right)
 $$
where $a$ is an $r \times r$ matrix, $b, d$ are $r \times t$ matrices, $c$ is an $s \times t$ matrix and $e$ is a $t \times t$ matrix.
Note that $a$ is invertible and $\left( \begin{array}{c} d \\ e \\ \end{array} \right)$ has full rank.
By construction, the representative of $x'$ is given
$$ x' = \left(
          \begin{array}{cc}
            0 & d \\
            0 & e \\
          \end{array}
        \right).
  $$
The matrix representative of $\overline{x}'$ may be written as
$ \overline{x}' = \left(
                     \begin{array}{cc}
                       f & g \\
                       0 & h \\
                     \end{array}
                   \right).
$
Note that $\ker {x}^2 / \ker x $ is invariant under $\overline{x}'$ by the assumption $[x', \overline{x}']=0$.
We have
\begin{align} \label{Eq: 0 = [x' olx']}
 0 = [x', \overline{x}'] = x' \overline{x}' - \overline{x}'x'  = \left(
                                                 \begin{array}{cc}
                                                   0 & dh - fd-ge \\
                                                   0 & eh - he \\
                                                 \end{array}
                                               \right).
\end{align}

Let $X := a f a^{-1}$ and $Y$ be an invertible $s \times s$ matrix. Since $\left( \begin{array}{c} d \\ e \\ \end{array} \right)$ has full rank,
the following equations
$$ Z_1 d + Z_3 e =  ag + bh -Xb \ \text{ and } \ Z_2 d + Z_4 e = ch-Yc  $$
have solutions. We choose a solution $ \left(\begin{array}{cc}
                                        Z_1 & Z_3 \\
                                        Z_2 & Z_4
                                      \end{array} \right)
$ of the above equations which maps $V_i$ to $V_{i-1}$ for $i\in I$.
We set $$\overline{x} := \left(
                      \begin{array}{cccc}
                        X & 0 & Z_1 & Z_3 \\
                        0 & Y & Z_2 & Z_4 \\
                        0 & 0 & f & g \\
                        0 & 0 & 0 & h \\
                      \end{array}
                    \right)
$$
By construction, we have $\overline{x}|_{M / \ker x} = \overline{x}' $ and
$$ [x, \overline{x}] = \left(
                         \begin{array}{cccc}
                           0 & 0 & af-Xa & ag + bh -Xb -Z_1 d - Z_3 e \\
                           0 & 0 & 0 & ch-Yc - Z_2 d - Z_4 e  \\
                           0 & 0 & 0 & 0 \\
                           0 & 0 & 0 & 0 \\
                         \end{array}
                       \right) = 0.
 $$

Now we take a map $ t \in \bigoplus_{i\in I} \Hom( M_i, L_i)$ such that, for any element $m = xm' \in \im\ x$,
$$ t(m) = t'(m'+\ker x).  $$
Then we have $ t\circ x |_{M/ \ker x} = t' $.
Fix an ordered basis of $L$, and let $\mathfrak{t}'$ be the matrix representative of $t'|_{\ker x^2/\ker x}$. Since
$\ker x' \cap \ker \overline{x}' \cap \ker t' = 0$, we have
$$ \ker f \cap \ker \mathfrak{t}' = 0. $$
Then the matrix representative $\mathfrak{t}$ of $t|_{\ker x}$ can be written as
$$ \mathfrak{t} = \left(
                    \begin{array}{cc}
                      \mathfrak{t}' a^{-1} & \mathfrak{t}'' \\
                    \end{array}
                  \right)
 $$
for some $r \times \ell$ matrix $\mathfrak{t}''$. Since $\overline{x}|_{\ker x} = \left(
                                                                                    \begin{array}{cc}
                                                                                      X & 0 \\
                                                                                      0 & Y \\
                                                                                    \end{array}
                                                                                  \right)
$, $X = a f a^{-1}$ and $Y$ is invertible, we have
$$ \ker (\overline{x}|_{\ker x}) \cap \ker (t|_{\ker x}) = 0, $$
which yields
$$ \ker x \cap \ker \overline{x} \cap \ker t = 0. $$
Therefore, we have the assertion.
\end{proof}


Recall the fundamental isomorphism of perfect crystals $\eqref{Eq:fundamental thm of perfect crystals}$. Combining
the perfect crystals $\Bone{\ell}$ and $\Ben{\ell}$ with Theorem \ref{Thm:saito geometric crystal}, we have
\begin{align*}
\Phi_\Lambda^1 &: \gB(\Lambda) \longrightarrow \gB(\Lambda')\otimes \Bone{\ell}, \\
\Phi_\Lambda^n &: \gB(\Lambda) \longrightarrow \gB(\Lambda'')\otimes \Ben{\ell},
\end{align*}
where $\Lambda' := \sum_{i\in I} \Lambda(h_{i+1})\Lambda_i$ and $\Lambda'' := \sum_{i\in I} \Lambda(h_{i-1})\Lambda_i$.
We give a geometric interpretation of the isomorphisms $\Phi_\Lambda^1$ and $\Phi_\Lambda^n$ in terms of quiver varieties.

\begin{Thm} \label{Thm:Fundamantal thm for qv}
Let $X \in \gB(\Lambda)$ and $\alpha := \Lambda -  \wt(X)$. Set
$$ \beta = \underline{\dim}(\ker x)\ \  (\text{resp.\ } \gamma = \underline{\dim}(\ker \overline{x}) ) $$
for $[\chi=x+\overline{x}, t]$ in a nonempty open subset of $X$ as in Lemma \ref{Lem:open set U in B1 and Bn}.
Let $\phi: V(\alpha)/ \ker x \to V(\alpha-\beta)$ (resp.\ $\bar{\phi}: V(\alpha)/ \ker \overline{x} \to V(\alpha-\gamma)$) be an $I$-graded vector space isomorphism.
\begin{enumerate}
\item Let $\Lambda' := \sum_{i\in I} \Lambda(h_{i+1})\Lambda_i$. Then there exists a unique irreducible component $X' \in \gB(\Lambda')$ satisfying the following conditions:
\begin{enumerate}
\item there exists a nonempty open subset $U \subset X$ such that
$$ [ \phi \circ (\chi|_{V(\alpha)/\ker x}) \circ \phi^{-1},\   (t \circ x|_{V(\alpha)/\ker x}) \circ \phi^{-1} ] \in X' $$
for any element $[\chi=x+\overline{x}, t] \in U$,
\item there exists a nonempty open subset $U' \subset X'$ such that every element $[\chi', t'] \in U'$ can be written as
$$ [\chi', t'] = [ \phi \circ (\chi|_{V(\alpha)/\ker x}) \circ \phi^{-1},  (t \circ x|_{V(\alpha)/\ker x})\circ \phi^{-1}]  $$
for some element $[\chi = x + \overline{x}, t]\in X$,
\item moreover, we have
$$ \Phi^1_\Lambda(X) = X' \otimes \psi^{1,\ell}(\wt(b^{\Lambda}_0) - \cl(\beta)), $$
where $\psi^{1,\ell}$ is the 1-1 map in $\eqref{Eq:map wt to B1}$,
and $\mathbf{p}^1(\Lambda) = ( b^{\Lambda}_k )_{k \ge 0}$ is the ground-state path of $\Lambda$ in $\Bone{\ell}$ given in $\eqref{Eq:g-s path in B1}$.
\end{enumerate}
\item Let $\Lambda'' := \sum_{i\in I} \Lambda(h_{i-1})\Lambda_i$. Then there exists a unique irreducible component $X'' \in \gB(\Lambda'')$ satisfying the following conditions:
\begin{enumerate}
\item there exists a nonempty open subset $U \subset X$ such that
$$ [ \bar{\phi} \circ (\chi|_{V(\alpha)/\ker \overline{x}}) \circ \bar{\phi}^{-1},\  (t \circ \overline{x}|_{V(\alpha)/\ker \overline{x}}) \circ \bar{\phi}^{-1}  ] \in X'', $$
for any element $[\chi=x+\overline{x}, t] \in U$,
\item there exists a nonempty open subset $U'' \subset X''$ such that every element $[\chi'', t''] \in U''$ can be written as
$$ [\chi'',t''] = [\bar{\phi} \circ (\chi|_{V(\alpha)/\ker \overline{x}}) \circ \bar{\phi}^{-1}, (t \circ \overline{x}|_{V(\alpha)/\ker \overline{x}}) \circ \bar{\phi}^{-1} ] $$
for some element $[\chi = x + \overline{x}, t]\in X$,
\item moreover, we have
$$ \Phi^n_\Lambda(X) = X'' \otimes \psi^{n,\ell}(\wt(\overline{b}^{\Lambda}_0) - \cl(\gamma)) , $$
where $\psi^{n,\ell}$ is the 1-1 map in $\eqref{Eq:map wt to Bn}$,
and $\mathbf{p}^n(\Lambda) = ( \overline{b}^{\Lambda}_k )_{k \ge 0}$ is the ground-state path of $\Lambda$ in $\Ben{\ell}$ given in $\eqref{Eq:g-s path in Bn}$.
\end{enumerate}
\end{enumerate}
\end{Thm}
Note that the open sets $U$, $U'$ and $U''$ in Theorem \ref{Thm:Fundamantal thm for qv} do not depend on the choices of $\phi$ and $\bar{\phi}$ since $[\chi,t]$ denotes the $G(\alpha)$-orbit of $(\chi, t)\in X(\Lambda, \alpha)$.
\begin{proof}
We first deal with the assertion (1).
Let $\mathbf{Y}=(Y^{(1)}, \ldots, Y^{(\ell)})$ be the element of $\mathcal{Y}^1(\Lambda)$ corresponding to $X$ via the isomorphism in Proposition \ref{Prop:isom between Y1 and gB}.
Recall that $(x(\mathbf{Y}), V(\alpha))$ is the representation of $(I, \Omega)$ appeared in $\eqref{Eq:def of x(Y)}$.

(a) From $\eqref{Eq:def of Pi}$ and $\eqref{Eq:def of irr component of Y}$, there is a nonempty open subset
$U \subset X$ such that, for any element $[\chi, t] \in U$, $\chi$ is contained in the conormal bundle of the $G(\alpha)$-orbit of $x(\mathbf{Y})$.
For $1 \le k \le \ell$, let $Y^{\prime(k) }$ be the Young wall obtained from $Y^{(k)}$ by removing the 0th column $y_0^{(k)}$ of $Y^{(k)}$.
It follows from Equation $\eqref{Eq:computation of ker Y}$ that
\begin{align}\label{Eq:gv of thm of pc 1}
\beta = \sum_{k=1}^\ell \wt(y_0^{(k)}).
\end{align}
By construction, $\mathbf{Y}' := (Y^{\prime(1)}, \ldots, Y^{\prime(\ell)})$ can be considered as the element in $\mathcal{Y}^1(\Lambda')$.
Set
$$X' := \text{the irreducible component in $\gB(\Lambda')$ corresponding to $\mathbf{Y}'$} .$$
Take an element $[\chi, t] \in U $, and let $\chi' := \chi|_{V(\alpha)/\ker x} $ and $t' := t \circ x|_{V(\alpha)/\ker x}$.
We write $\chi' = x' + \overline{x}' $.  Note that $\chi'$ and $t'$ are well-defined by Lemma \ref{Lem:chi statble}.
By construction, $\phi \circ x' \circ \phi^{-1}$ is contained in the $G(\alpha - \beta)$-orbit of $x(\mathbf{Y}')$.
Since $ x'\overline{x}' -  \overline{x}'x' = (x \overline{x} - \overline{x} x)|_{V(\alpha)/\ker x} = 0$, by $\eqref{Eq:commuting property of chi}$, we have
$$\mu_i (\phi \circ \chi' \circ \phi^{-1}) = 0 \ (i\in I).$$
Hence $[\phi \circ \chi' \circ \phi^{-1},\  t' \circ \phi^{-1}  ]$ is contained in $X'$ if $(\chi', t')$ satisfies the stability condition.

Suppose that there is a nonzero $\chi'$-invariant subspace $S'$ of $V(\alpha)/\ker x$ with $t' (S')=0$.
Let $\pi: V(\alpha) \to V(\alpha)/\ker x$ be the natural projection. Set
$$S:= x( \pi^{-1}(S')) \subset V(\alpha). $$
Since $S'$ is nonzero, $S$ is a nonzero $I$-graded subspace of $V(\alpha)$. For any element $v \in S$,
$v$ can be written as $v = x w$ for some $w \in \pi^{-1}(S')$. Note that $\chi w \in \pi^{-1}(S')$ since $S'$ is $\chi'$-invariant. Then, it follows from
$$ \chi v = (x + \overline{x})x w = x \chi(w)$$
that $S$ is $\chi$-invariant. Moreover, by the assumption $t' (S')=0$, we have $t(S)=0$. Hence
$S$ is a nonzero $\chi$-invariant subspace with $t(S)=0$, which is a contradiction. Therefore,
$(\chi', t')$ satisfies the stability condition.

(b) Let $U'$ be a nonempty open subset of $X'$ such that, for any element $[\chi', t'] \in U'$, $\chi'$ is contained in the conormal bundle of the $G(\alpha-\beta)$-orbit of $x(\mathbf{Y}')$.
Take an element $[\chi'=x' + \overline{x}', t'] \in U'$. By construction,
$x'$ is contained in the $G(\alpha- \beta)$-orbit of $x(\mathbf{Y}')$. Hence,
there exists an element $x$ in the $G(\alpha)$-orbit of $x(\mathbf{Y})$ such that
$$ x' = \phi \circ x|_{V(\alpha)/\ker x} \circ \phi^{-1}. $$
Then the assertion (b) follows from Lemma \ref{Lem:existance of the open set}.

(c) Combining Proposition \ref{Prop:iso of YW and Path} with the fundamental isomorphism of perfect crystals $\eqref{Eq:fundamental thm of perfect crystals}$,
we have
$$ \Phi: \mathcal{Y}^1(\Lambda) \longrightarrow \mathcal{Y}^1(\Lambda') \otimes \Bone{\ell}. $$
Then, by Equation $\eqref{Eq:gv of thm of pc 1}$ and the isomorphism $\eqref{Eq:iso between Y1 and P1}$, we obtained
$$ \Phi(\mathbf{Y}) = \mathbf{Y}' \otimes \psi^{1,\ell}( \wt(b_0^\Lambda)) - \cl(\beta)), $$
which yields
$$ \Phi^1_\Lambda(X) = X' \otimes b .  $$

The remaining case (2) can be obtained in the same manner.
\end{proof}

Let $\alpha, \beta \in Q^+$ with $\beta \le \alpha$, and $\Lambda' = \sum_{i\in I}\Lambda(h_{i+1})\Lambda_i$.
We consider an analogue of the diagram given in \cite[Section 12.10]{L91}:
\begin{align} \label{Eq:diagram of X(L,a)}
X(\Lambda', \alpha-\beta) \buildrel p_1 \over \longleftarrow F' \buildrel p_2 \over \longrightarrow F''  \buildrel p_3 \over \longrightarrow X(\Lambda, \alpha),
\end{align}
where $F''$ is the variety of all triples $(\chi=x+\overline{x}, t, M)$ such that
\begin{enumerate}
\item[(a)] $(\chi, t) \in X(\Lambda, \alpha)$,
\item[(b)] $M$ is a $\chi$-invariant subspace of $V(\alpha)$ with $\underline{\dim}M=\beta$ and $t\circ x(M)=0$,
\end{enumerate}
and $F'$ is the variety of all quintuples $(\chi, t, M, f, g)$ such that
\begin{enumerate}
\item[(a)] $(\chi, t, M) \in F''$,
\item[(b)] $f=(f_i)_{i\in I}, g=(g_i)_{i\in I}$ give an exact sequence
$$ 0 \longrightarrow V_i(\beta) \buildrel f_i \over \longrightarrow V_i(\alpha) \buildrel g_i \over \longrightarrow V_i(\alpha-\beta) \longrightarrow 0  $$
such that $\im\ f = M$.
\end{enumerate}
We also have
$$ p_1( \chi, t, M, f, g) := ( \tilde{g} \circ \chi|_{V(\alpha)/M} \circ  \tilde{g}^{-1}, \  (t\circ x|_{V(\alpha)/M}) \circ \tilde{g}^{-1}  ), $$
where $\chi = x+\overline{x}$, and $\tilde{g}:V(\alpha)/M \rightarrow V(\alpha-\beta)$ is the $I$-graded isomorphism induced by $g$,
$$ p_2(\chi, t, M, f,g) := (\chi, t, M)\ \text{ and }\ p_3(\chi, t, M):= (\chi, t). $$
Note that $p_2$ is a $G(\alpha-\beta)\times G(\beta)$-principle bundle and an open map.

Take a nonempty open subset $U$ of $X \in \gB(\Lambda)$ as in Lemma \ref{Lem:open set U in B1 and Bn} and let $\alpha:= \Lambda - \wt(X)$
and $\beta := \underline{\dim}(\ker x)$ for $[\chi = x+\overline{x},t] \in U$.
Let
$$\pi: X(\Lambda, \alpha)^{\rm st} \rightarrow X(\Lambda, \alpha)^{\rm st}/G(\alpha)$$
and $U_0 := \pi^{-1}(U)$. Note that $X(\Lambda, \alpha)^{\rm st} $ is
a nonempty open subset of $X(\Lambda, \alpha)$. Define a map $\iota:U_0  \rightarrow F'' $ by
\begin{align*}
 \iota(\chi, t) = (\chi, t, \ker x)
\end{align*}
for $(\chi=x+\overline{x}, t) \in U_0$. Note that $p_3 \circ \iota = {\rm id}|_{U_0}$. We set $X'\in \gB(\Lambda')$ to be
the irreducible component given in Theorem \ref{Thm:Fundamantal thm for qv} (1). Given a nonempty open subset $U_0' \subset X(\Lambda', \alpha-\beta)$
with $U_0' \cap \pi^{-1}(X') \ne \emptyset$, it follows from Lemma \ref{Lem:existance of the open set} that
$$ \tilde{U}_0 := \iota^{-1} \circ p_2 \circ p_1^{-1}(U_0') $$
is a nonempty open subset of $U_0$. Therefore, given a nonempty open subset $U' \subset X'$, there exists a nonempty open subset $\tilde{U} \subset X$ such that,
for any element $[\chi=x+ \overline{x}, t] \in \tilde{U}$,
$$ [\phi \circ (\chi|_{V(\alpha)/\ker x}) \circ \phi^{-1}, (t\circ x |_{V(\alpha)/\ker x}) \circ \phi^{-1}] \in U' ,$$
where $\phi:V(\alpha)/\ker x \rightarrow V(\alpha-\beta)$ is an $I$-graded vector space isomorphism.

In the same manner, let $\Lambda'' = \sum_{i\in I}\Lambda(h_{i-1})\Lambda_i$ and $\gamma := \underline{\dim}(\ker \overline{x})$ for $[\chi = x+\overline{x},t] \in U$. We consider the following diagram:
$$X(\Lambda'', \alpha-\gamma) \buildrel \overline{p}_1 \over \longleftarrow \overline{F}' \buildrel \overline{p}_2 \over \longrightarrow \overline{F}''  \buildrel \overline{p}_3 \over \longrightarrow X(\Lambda, \alpha),$$
where $\overline{F}''$ is the variety of all triples $(\chi=x+\overline{x}, t, M)$ such that
\begin{enumerate}
\item[(a)] $(\chi, t) \in X(\Lambda, \alpha)$,
\item[(b)] $M$ is a $\chi$-invariant subspace of $V(\alpha)$ with $\underline{\dim}M=\gamma$ and $t\circ \overline{x}(M)=0$,
\end{enumerate}
and $\overline{F}'$ is the variety of all quintuples $(\chi, t, M, f, g)$ such that
\begin{enumerate}
\item[(a)] $(\chi, t, M) \in \overline{F}''$,
\item[(b)] $f=(f_i)_{i\in I}, g=(g_i)_{i\in I}$ give an exact sequence
$$ 0 \longrightarrow V_i(\gamma) \buildrel f_i \over \longrightarrow V_i(\alpha) \buildrel g_i \over \longrightarrow V_i(\alpha-\gamma) \longrightarrow 0  $$
such that $\im\ f = M$.
\end{enumerate}
We set
$$ \overline{p}_1( \chi, t, M, f, g) := ( \tilde{g} \circ \chi|_{V(\alpha)/M} \circ  \tilde{g}^{-1}, \  (t\circ \overline{x}|_{V(\alpha)/M}) \circ \tilde{g}^{-1}  ), $$
where $\chi = x+\overline{x}$, and $\tilde{g}:V(\alpha)/M \rightarrow V(\alpha-\gamma)$ is the $I$-graded isomorphism induced by $g$,
$$ \overline{p}_2(\chi, t, M, f,g) := (\chi, t, M)\ \text{ and }\ \overline{p}_3(\chi, t, M):= (\chi, t). $$
Define a map $ \overline{\iota}:U \rightarrow \overline{F}'' $ by
$$ \overline{\iota}(\chi, t)=(\chi, t, \ker \overline{x})  $$
for $[\chi=x+\overline{x}, t] \in U$, and let $X''$ be the irreducible component given in Theorem \ref{Thm:Fundamantal thm for qv} (2). Then we obtain that, given
a nonempty open subset $U'' \subset X''$, there exists a nonempty open subset $\tilde{U} \subset X$ such that
for any element $[\chi=x+ \overline{x}, t] \in \tilde{U}$,
$$ [\bar{\phi} \circ (\chi|_{V(\alpha)/\ker \overline{x}}) \circ \bar{\phi}^{-1}, (t\circ \overline{x} |_{V(\alpha)/\ker \overline{x}}) \circ \bar{\phi}^{-1}] \in U'', $$
where $\bar{\phi}:V(\alpha)/\ker \overline{x} \rightarrow V(\alpha-\gamma)$ is an $I$-graded vector space isomorphism .
Consequently, we obtain the following lemma.

\begin{Lem} \label{Lem:open to open}
With the same notations as in Theorem \ref{Thm:Fundamantal thm for qv}, we have the following.
\begin{enumerate}
\item  Given a nonempty open subset $U' \subset X'$, there exists a nonempty open subset $\tilde{U} \subset X$ such that
$$ [\phi \circ (\chi|_{V(\alpha)/\ker x}) \circ \phi^{-1}, (t\circ x |_{V(\alpha)/\ker x}) \circ \phi^{-1}] \in U' $$
for any element $[\chi=x+ \overline{x}, t] \in \tilde{U}$.
\item Given a nonempty open subset $U'' \subset X''$, there exists a nonempty open subset $\tilde{U} \subset X$ such that
$$ [\bar{\phi} \circ (\chi|_{V(\alpha)/\ker \overline{x}}) \circ \bar{\phi}^{-1}, (t\circ \overline{x} |_{V(\alpha)/\ker \overline{x}}) \circ \bar{\phi}^{-1}] \in U'' $$
for any element $[\chi=x+ \overline{x}, t] \in \tilde{U}$.
\end{enumerate}
\end{Lem}

Let us recall the path realizations $\mathcal{P}^1(\Lambda)$, $\mathcal{P}^n(\Lambda)$ and $\mathcal{P}^{\rm ad}(\Lambda)$
arising from the perfect crystals $\Bone{\ell}, \Ben{\ell}$ and $\Bad{\ell}$ respectively. By Theorem \ref{Thm: crystal iso of path realization} and Theorem \ref{Thm:saito geometric crystal},
we have crystal isomorphisms from $\gB(\Lambda)$ to $\mathcal{P}^1(\Lambda)$, $\mathcal{P}^n(\Lambda)$ and $\mathcal{P}^{\rm ad}(\Lambda)$ respectively.
We give explicit descriptions in terms of quiver varieties.

\begin{Thm} \label{Thm:iso from gB to P1 and Pn}

Let
$$ \Upsilon^1:\gB(\Lambda) \longrightarrow \mathcal{P}^1(\Lambda) \ \  (resp.\ \Upsilon^n:\gB(\Lambda) \longrightarrow \mathcal{P}^n(\Lambda) ) $$
be the unique crystal isomorphism given by Theorem \ref{Thm: crystal iso of path realization} and Theorem \ref{Thm:saito geometric crystal}.
Let $\mathbf{p}^1(\Lambda) := ( b^{\Lambda}_i )_{i \ge 0}$ (resp.\ $\mathbf{p}^n(\Lambda) := ( \overline{b}^{\Lambda}_i )_{i \ge 0}$)
be the ground-state path of $\Lambda$ in $\Bone{\ell}$ (resp.\ $\Ben{\ell}$) given in $\eqref{Eq:g-s path in B1}$ (resp.\ $\eqref{Eq:g-s path in Bn}$).
Take an irreducible component $X \in \gB(\Lambda)$. Then there exists a nonempty open subset $U \subset X$ satisfying the following conditions:
\begin{enumerate}
\item  let
\begin{align*}
\mathfrak{r}_i &:= \wt(b_i^{\Lambda}) - \cl(\underline{\dim} (\ker x^{i+1}/\ker x^i) ), \\
\mathfrak{s}_i &:= \wt(\overline{b}_i^{\Lambda}) - \cl(\underline{\dim} (\ker \overline{x}^{i+1}/\ker \overline{x}^i) )
\end{align*}
for $[\chi=x+\overline{x}, t] \in U$, then the weights $\mathfrak{r}_i$ and $\mathfrak{s}_i$ do not depend on the choice of $[\chi, t] $,

\item moreover,
$$ \Upsilon^1(X) = ( \psi^{1,\ell}(\mathfrak{r}_i) )_{i\ge 0 } \ \text{ and }\  \Upsilon^n(X) = ( \psi^{n,\ell}(\mathfrak{s}_i) )_{i\ge 0 }, $$
where $\psi^{1,\ell}$ (resp.\ $\psi^{n,\ell}$) is the 1-1 map in $\eqref{Eq:map wt to B1}$ (resp.\ $\eqref{Eq:map wt to Bn}$).

\end{enumerate}

\end{Thm}

\begin{proof}

Let $\alpha := \Lambda - \wt(X)$. We will use induction on $\HT(\alpha)$. If $\HT(\alpha) = 0$, then there is nothing to prove.
Assume that $\HT(\alpha) > 0$. We first deal with the case of $\Upsilon^1$ and $\ker x$. Let $\Lambda' := \sum_{i\in I} \Lambda(h_{i+1})\Lambda_i$,
and $X'$ be the irreducible component of $\gB(\Lambda')$ given in Theorem \ref{Thm:Fundamantal thm for qv} (1).
By the induction hypothesis, there is a nonempty open subset $U' \subset X'$ satisfying the conditions of our assertion. By Lemma \ref{Lem:open to open}, there is a nonempty open subset
$U \subset X$ such that
$$ \underline{\dim} \ker x'^{i+1}/ \ker x'^{i} = \underline{\dim} \ker x^{i+2}/ \ker x^{i+1}   $$
for $i \in \Z_{> 0 }$, $[\chi'=x' + \overline{x}', t'] \in U'$ and $[\chi=x + \overline{x}, t] \in U$.
Since the subsequence $( b^{\Lambda}_i )_{i \ge 1}$ of the ground-state path $\mathbf{p}^1(\Lambda) $ can be viewed as the ground-state path of $\Lambda'$ in $\Bone{\ell}$,
 it follows from Theorem \ref{Thm:Fundamantal thm for qv} (1) that $\mathfrak{r}_i$ does not depend on the choice of $[\chi, t] \in U$ for each $i$, and
$$\Upsilon^1(X) = ( \psi^{1,\ell}(\mathfrak{r}_i) )_{i\ge 0 }.$$

The remaining case, $\Upsilon^n$ and $\ker \overline{x}$, can be proved in the same manner.
\end{proof}

\begin{Thm}  \label{Thm:iso from gB to Pad}
Let
$$ \Upsilon^{\rm ad}:\gB(\Lambda) \longrightarrow \mathcal{P}^{\rm ad}(\Lambda)  $$
be the unique crystal isomorphism given by Theorem \ref{Thm: crystal iso of path realization} and Theorem \ref{Thm:saito geometric crystal}.
Let $\Lambda'' = \sum_{i\in I}\Lambda(h_{i-1})\Lambda_i$, and let $\mathbf{p}^1(\Lambda'') := ( b^{\Lambda''}_i )_{i \ge 0}$ (resp.\ $\mathbf{p}^n(\Lambda) := ( \overline{b}^{\Lambda}_i )_{i \ge 0}$)
be the ground-state path of $\Lambda''$ (resp.\ $\Lambda$) in $\Bone{\ell}$ (resp.\ $\Ben{\ell}$) given in $\eqref{Eq:g-s path in B1}$ (resp.\ $\eqref{Eq:g-s path in Bn}$).
Take an irreducible component $X \in \gB(\Lambda)$. Then there exists a nonempty open subset $U \subset X$ satisfying the following conditions:
\begin{enumerate}
\item let
\begin{align*}
\mathfrak{r}_i  &:= \wt(b_0^{\Lambda''}) - \cl(\underline{\dim} (\ker (x\overline{x})^{i+1}/\ker \overline{x}(x\overline{x})^i) ), \\
\mathfrak{s}_i &:= \wt(\overline{b}_0^{\Lambda}) - \cl(\underline{\dim} (\ker \overline{x}(x\overline{x})^{i}/\ker (x\overline{x})^i) )
\end{align*}
for $[\chi=x+\overline{x}, t] \in U$, then the weights $\mathfrak{r}_i$ and $\mathfrak{s}_i$ do not depend on the choice of $[\chi, t] $,

\item moreover,
$$ \Upsilon^{\rm ad}(X) = ( \psi^{{\rm ad},\ell}(\mathfrak{r}_i, \mathfrak{s}_i) )_{i\ge 0 }, $$
where $\psi^{{\rm ad},\ell}$ is the 1-1 map in $\eqref{Eq:map wt to Bad}$.

\end{enumerate}
\end{Thm}

\begin{proof}
Let $\alpha := \Lambda - \wt(X)$. We will use induction on $\HT(\alpha)$. Since the case $\HT(\alpha) = 0$ is trivial,
we may assume $\HT(\alpha) > 0$. Let $\gamma := \underline{\dim} (\ker \overline{x}) $ for $[\chi=x+\overline{x}, t]$
in a nonempty open subset of $X$ as in Lemma \ref{Lem:open set U in B1 and Bn}, and let $X''$ be the irreducible component of $\gB(\Lambda'')$ associated
with $X$ as in Theorem \ref{Thm:Fundamantal thm for qv} (2). Similarly, let $\beta := \underline{\dim} (\ker x'') $ for $[\chi''=x''+\overline{x}'', t'']$
in a nonempty open subset of $X''$ as in Lemma \ref{Lem:open set U in B1 and Bn}, and let $X'$ be the irreducible component of $\gB(\Lambda)$ associated
with $X''$ as in Theorem \ref{Thm:Fundamantal thm for qv} (1). Then we have
$$ \Phi^1_{\Lambda ''}(X'') = X' \otimes b \ \text{ and } \  \Phi^n_{\Lambda }(X) = X'' \otimes \overline{b} $$
for some $b \in \Bone{\ell}$ and $\overline{b} \in \Ben{\ell}$.
From the fundamental isomorphism $\eqref{Eq:fundamental thm of perfect crystals}$ and Lemma \ref{Lem:g-s path in Bad}, there is a crystal isomorphism
$$ \Phi^{\rm ad}_\Lambda: \gB(\Lambda) \rightarrow \gB(\Lambda) \otimes \Bad{\ell}. $$
By Lemma \ref{Lem:B1 Bn isom Bad}, we obtain
$$ \Phi^{\rm ad}_\Lambda(X) = X' \otimes \xi^{\ell}(b \otimes \overline{b}). $$

On the other hand, by the induction hypothesis, there is an open subset $U' \subset X'$ satisfying the conditions of our assertion. Then it follows
from Lemma \ref{Lem:open to open} that there is a nonempty open subset $U \subset X$ such that, for $[\chi = x+\overline{x}, t] \in U$,
$$ [ \tilde{\phi} \circ ( \chi |_{V(\alpha)/\ker x \overline{x}}) \circ \tilde{\phi}^{-1}, (t\circ x  \overline{x}|_{V(\alpha)/\ker x \overline{x}}) \circ \tilde{\phi}^{-1} ] \in U', $$
where $ \tilde{\phi}: V(\alpha)/\ker x \overline{x} \rightarrow V(\alpha-\beta-\gamma)  $ is an $I$-graded vector space.
Note that $\gamma = \underline{\dim} (\ker \overline{x})$, $\beta = \underline{\dim}(\ker x\overline{x}/\ker \overline{x})$ and $\beta + \gamma = \underline{\dim}(\ker x\overline{x})$.
 By Theorem \ref{Thm:Fundamantal thm for qv}, we have
$$\wt(\overline{b}) = \wt(b_0^{\Lambda}) - \cl(\gamma)\ \text{ and } \ \wt(b) = \wt(b_0^{\Lambda''}) - \cl(\beta), $$
which yields, by Equation $\eqref{Eq:map wt to Bad}$,
$$ \xi^{\ell}(b \otimes \overline{b}) = \psi^{{\rm ad},\ell}(\mathfrak{r}_0, \mathfrak{s}_0) .$$
Since
\begin{align*}
\underline{\dim} (\ker (x\overline{x})^{i+2}/\ker \overline{x}(x\overline{x})^{i+1}) = \underline{\dim} (\ker (x'\overline{x}')^{i+1}/\ker \overline{x}'(x'\overline{x}')^i), \\
\underline{\dim} (\ker \overline{x}(x\overline{x})^{i+1}/\ker \overline{x}(x\overline{x})^{i+1}) = \underline{\dim} \overline{x}'(\ker (x'\overline{x}')^{i}/\ker \overline{x}'(x'\overline{x}')^i) \\
\end{align*}
for $[\chi=x+\overline{x},t] \in U$ and $[\chi'=x'+\overline{x}',t'] \in U'$,
our assertion follows from a standard induction argument.
\end{proof}

The following corollary, which is a consequence of Theorem \ref{Thm:Fundamantal thm for qv} and Theorem \ref{Thm:iso from gB to Pad}, can be regarded
as a geometric interpretation of the fundamental isomorphism of perfect crystals
$$ \Phi_\Lambda^{\rm ad}: \gB(\Lambda) \buildrel \sim \over \longrightarrow \gB(\Lambda) \otimes \Bad{\ell}. $$

\begin{Cor} \label{Cor: fundamental thm for ad}
Let $X \in \gB(\Lambda)$ and $\alpha := \Lambda - \wt(X)$. Set
$$ \beta :=   \underline{\dim} \ker x \overline{x} / \ker \overline{x} \ \text{ and }\
   \gamma :=   \underline{\dim} \ker \overline{x} $$
for $[\chi=x+\overline{x}]$ in a nonempty open subset of $X$ as in Theorem \ref{Thm:iso from gB to Pad} (1).
Let $ \tilde{\phi}: V(\alpha)/\ker x \overline{x} \rightarrow V(\alpha-\beta-\gamma)$ be an $I$-graded vector space.
Then there is a unique irreducible component $X' \in \gB(\Lambda)$ satisfying the following conditions:
\begin{enumerate}
\item there is a nonempty open subset $U \subset X$ such that
$$ [\tilde{\phi} \circ (\chi|_{V(\alpha)/\ker x\overline{x}})\circ \tilde{\phi}^{-1}, (t \circ x\overline{x}|_{V(\alpha)/ \ker x\overline{x}} )\circ \tilde{\phi}^{-1}] \in X' $$
for any element $[\chi=x+\overline{x}, t] \in U$,
\item there is a nonempty open subset $U' \subset X'$ such that every element $[\chi', t']\in U'$ can be written as
$$ [\chi', t'] = [\tilde{\phi} \circ (\chi|_{V(\alpha)/\ker x\overline{x}})\circ \tilde{\phi}^{-1}, (t \circ x\overline{x}|_{V(\alpha)/\ker x\overline{x}})\circ \tilde{\phi}^{-1}]$$
for some element $[\chi=x+\overline{x},t]\in X$,
\item moreover,
$$ \Phi_\Lambda^{\rm ad}(X) = X' \otimes \psi^{{\rm ad}, \ell}(\wt(b_0^{\Lambda''})-\cl(\beta), \wt(\overline{b}_0^{\Lambda})-\cl(\gamma)), $$
where $\psi^{{\rm ad}, \ell}$ is the 1-1 map in $\eqref{Eq:map wt to Bad}$, and $\mathbf{p}^1(\Lambda'') := ( b^{\Lambda''}_i )_{i \ge 0}$ (resp.\ $\mathbf{p}^n(\Lambda) := ( \overline{b}^{\Lambda}_i )_{i \ge 0}$)
is the ground-state path of $\Lambda''$ (resp.\ $\Lambda$) in $\Bone{\ell}$ (resp.\ $\Ben{\ell}$) given in $\eqref{Eq:g-s path in B1}$ (resp.\ $\eqref{Eq:g-s path in Bn}$).
\end{enumerate}
\end{Cor}
Note that the open sets $U$ and $U'$ in Corollary \ref{Cor: fundamental thm for ad} do not depend on the choice of $\tilde{\phi}$ since $[\chi,t]$ denotes the $G(\alpha)$-orbit of $(\chi, t)\in X(\Lambda, \alpha)$.

\vskip 1em

\begin{Exa}
We keep all notations in Example \ref{Ex:Young walls}.
 By $\eqref{Eq:def of Pi}$ and $\eqref{Eq:def of irr component of Y}$,
it suffices to consider the fiber $\pi_\Omega^{-1}(x(\mathbf{Y}))$. By \cite[Section 12.8, Proposition 15.5]{L91}, we have
\begin{align*}
\pi^{-1}_\Omega(x(\mathbf{Y})) &= \{x(\mathbf{Y}) + \overline{x} \mid \overline{x} \in E_{\overline{\Omega}}(\alpha),\ [x(\mathbf{Y}), \overline{x}]=0 \} \\
&= \{x(\mathbf{Y}) + \overline{x} \mid \overline{x} = \overline{x}( a_1, \cdots, a_{29}),\ a_1, \ldots, a_{29} \in \C \}.
\end{align*}
Here,
\begin{align*}
\overline{x}( a_1, \cdots, a_{29} ) := \left(
                                                          \begin{array}{ccc}
                                                            0 & \overline{x}_1 & 0  \\
                                                            0 & 0 & \overline{x}_2  \\
                                                            \overline{x}_0 & 0 & 0  \\
                                                          \end{array}
                                                        \right) \in E_{\overline{\Omega}}(\alpha),
\end{align*}
\begin{align*}
\overline{x}_0 &:=\left(
                   \begin{array}{cccc}
                     0 & 0 & 0 & 0 \\
                     0 & 0 & 0 & 0 \\
                     0 & 0 & 0 & 0 \\
                     0 & 0 & 0 & 0 \\
                     0 & 0 & 0 & a_{25} \\
                     0 & 0 & 0 & 0 \\
                   \end{array}
                 \right),\
\
\overline{x}_1 := \left(
                    \begin{array}{ccccccc}
                      a_1 & a_2 & a_3 & a_4 & 0 & a_5 & a_6 \\
                      a_7 & a_8 & a_9 & a_{10} & 0 & a_{11} & a_{12} \\
                      a_{26} & 0 & a_{27} & 0 & 0 & a_{28} & a_{29} \\
                      a_{13} & a_{14} & a_{15} & a_{16} & 0 & a_{17} & a_{18} \\
                    \end{array}
                  \right),\\
\overline{x}_2 &:= \left(
                    \begin{array}{cccccc}
                      0 & 0 & 0 & 0 & 0 & 0 \\
                      0 & 0 & a_{19} & 0 & 0 & a_{20} \\
                      0 & 0 & 0 & 0 & 0 & 0 \\
                      0 & 0 & a_{21} & 0 & 0 & a_{22} \\
                      a_{26} & a_{27} & a_{23} & a_{28} & a_{29} & a_{24} \\
                      0 & 0 & 0 & 0 & 0 & 0 \\
                      0 & 0 & 0 & 0 & 0 & a_{25} \\
                    \end{array}
                  \right)
\end{align*}
for $a_1, \ldots, a_{29} \in \C$. Let
$$x:= x(\mathbf{Y}), \quad \overline{x} := \overline{x}( a_1, \cdots, a_{29}) $$
and consider $a_1, \ldots, a_{29}$ as indeterminates. Then we have
$$ \underline{\dim}(\ker x^k) = \left\{
                                  \begin{array}{ll}
                                    0 & \hbox{ if } k=0, \\
                                    3\alpha_0 + 3 \alpha_1 + 2 \alpha_2 & \hbox{ if } k=1, \\
                                    4\alpha_0 + 4 \alpha_1 + 5 \alpha_2 & \hbox{ if } k=2, \\
                                    4\alpha_0 + 6 \alpha_1 + 6 \alpha_2 & \hbox{ if } k=3, \\
                                    4\alpha_0 + 7\alpha_1 + 6\alpha_2 & \hbox{ otherwise,}
                                  \end{array}
                                \right.
 \text{ and } \
\underline{\dim}(\ker \overline{x}^k) = \left\{
                                          \begin{array}{ll}
                                            0 & \hbox{ if } k=0, \\
                                            3\alpha_0 + 3\alpha_1 + 3\alpha_2 & \hbox{ if } k=1, \\
                                            3\alpha_0 + 6\alpha_1 + 4\alpha_2 & \hbox{ if } k=2, \\
                                            4\alpha_0 + 6\alpha_1 + 5\alpha_2 & \hbox{ if } k=3, \\
                                            4\alpha_0 + 7\alpha_1 + 5\alpha_2 & \hbox{ if } k=4, \\
                                            4\alpha_0 + 7\alpha_1 + 6\alpha_2 & \hbox{ otherwise.}
                                          \end{array}
                                        \right.
$$
Hence, by Theorem \ref{Thm:iso from gB to P1 and Pn}, we obtain
\begin{align*}
\Upsilon^1(X) & = \cdots \otimes
\begin{tabular}{|c|c|c|} \hline 2 & 2 & 3 \\ \hline \end{tabular} \otimes
\begin{tabular}{|c|c|c|} \hline 2 & 3 & 3 \\ \hline \end{tabular} \otimes
\begin{tabular}{|c|c|c|} \hline 2 & 2 & 3 \\ \hline \end{tabular} \otimes
\begin{tabular}{|c|c|c|} \hline 3 & 3 & 3 \\ \hline \end{tabular} \otimes
\begin{tabular}{|c|c|c|} \hline 1 & 2 & 3 \\ \hline \end{tabular}\ , \\
\Upsilon^n(X) & = \cdots \otimes
\begin{tabular}{|c|c|c|} \hline $\overline{3}$ & $\overline{3}$ & $\overline{1}$ \\ \hline \end{tabular} \otimes
\begin{tabular}{|c|c|c|} \hline $\overline{2}$ & $\overline{2}$ & $\overline{2}$ \\ \hline \end{tabular} \otimes
\begin{tabular}{|c|c|c|} \hline $\overline{1}$ & $\overline{1}$ & $\overline{1}$ \\ \hline \end{tabular} \otimes
\begin{tabular}{|c|c|c|} \hline $\overline{3}$ & $\overline{3}$ & $\overline{2}$ \\ \hline \end{tabular} \otimes
\begin{tabular}{|c|c|c|} \hline $\overline{1}$ & $\overline{1}$ & $\overline{1}$ \\ \hline \end{tabular}\ \otimes
\begin{tabular}{|c|c|c|} \hline $\overline{2}$ & $\overline{1}$ & $\overline{1}$ \\ \hline \end{tabular}\ .
\end{align*}

On the other hand, the ground-state path of $\Lambda'' = 2\Lambda_1 + \Lambda_2$ in $\Bone{\ell}$ is given as follows.
$$ \cdots \otimes  \begin{tabular}{|c|c|c|} \hline 1 & 1 & 2 \\ \hline \end{tabular} \otimes
\begin{tabular}{|c|c|c|} \hline 2 & 2 & 3 \\ \hline \end{tabular} \otimes
\begin{tabular}{|c|c|c|} \hline 1 & 3 & 3 \\ \hline \end{tabular} \otimes
\begin{tabular}{|c|c|c|} \hline 1 & 1 & 2 \\ \hline \end{tabular}\ .$$
By a direct computation, we have
\begin{align*}
\underline{\dim}(\ker (x\overline{x})^k) &= \left\{
                                  \begin{array}{ll}
                                    0 & \hbox{ if } k=0, \\
                                    4\alpha_0 + 6 \alpha_1 + 5 \alpha_2 & \hbox{ if } k=1, \\
                                    4\alpha_0 + 7\alpha_1 + 6\alpha_2 & \hbox{ otherwise,}
                                  \end{array}
                                \right. \\
\underline{\dim}(\ker \overline{x}(x\overline{x})^k) &= \left\{
                                          \begin{array}{ll}
                                            3\alpha_0 + 3\alpha_1 + 3\alpha_2 & \hbox{ if } k=0, \\
                                            4\alpha_0 + 7\alpha_1 + 5\alpha_2 & \hbox{ if } k=1, \\
                                            4\alpha_0 + 7\alpha_1 + 6\alpha_2 & \hbox{ otherwise.}
                                          \end{array}
                                        \right.
\end{align*}
By Theorem \ref{Thm:iso from gB to Pad}, we obtain
\begin{align*}
\Upsilon^{\rm ad}(X)  = \cdots \otimes
\begin{tabular}{|c|c|} \hline 1 & 2 \\ \cline{1-2} 3  \\ \cline{1-1} \end{tabular} \otimes
\begin{tabular}{|c|c|} \hline 1 & 2 \\ \cline{1-2} 3  \\ \cline{1-1} \end{tabular} \otimes
\begin{tabular}{|c|c|} \hline 2 & 3 \\ \cline{1-2} 3  \\ \cline{1-1} \end{tabular}\otimes
\begin{tabular}{|c|c|c|c|c|c|}
  \hline
  1 & 2 & 2 & 2 & 2 & 3 \\
 \hline
  3 & 3 & 3  \\
  \cline{1-3}
\end{tabular} \ .
\end{align*}

\end{Exa}

\vskip 2em


\bibliographystyle{amsplain}

\begin{thebibliography}{10}

\bibitem{BFKL06}
G.~Benkart, I.~Frenkel, S.-J. Kang, and H.~Lee, \emph{Level 1 perfect crystals
  and path realizations of basic representations at $q=0$}, Int. Math. Res.
  Not. (2006), 1--28.

\bibitem{DJKMO89}
E.~Date, M.~Jimbo, A.~Kuniba, T.~Miwa, and M.~Okado, \emph{Paths, {M}aya
  diagrams and representations of $\widehat{{sl}}(r,\mathbb{C})$}, Adv. Stud.
  Pure Math. \textbf{19} (1989), 149--191.

\bibitem{FLOTW99}
O.~Foda, B.~Leclerc, M.~Okado, J.-Y. Thibon, and T.~A. Welsh, \emph{Branching
  functions of ${A}^{(1)}_{n-1}$ and {J}antzen-{S}eitz problem for
  {A}riki-{K}oike algebras}, Adv. Math. \textbf{141} (1999), no.~2, 322--365.

\bibitem{FS03}
I.~B. Frenkel and A.~Savage, \emph{Bases of representations of type ${A}$
  affine {L}ie algebras via quiver varieties and statistical mechanics}, Int.
  Math. Res. Not. (2003), no.~28, 1521--1547.

\bibitem{HK02}
J.~Hong and S.-J. Kang, \emph{{I}ntroduction to {Q}uantum {G}roups and
  {C}rystal {B}ases}, Grad. Stud. Math., 42, American Mathematical Society,
  Providence, RI, 2002.

\bibitem{JL09}
N.~Jacon and C.~Lecouvey, \emph{{K}ashiwara and {Z}elevinsky involutions in
  affine type {A}}, Pacific J. Math. \textbf{243} (2009), no.~2, 287--311.

\bibitem{JMMO91}
M.~Jimbo, K.~C. Misra, T.~Miwa, and M.~Okado, \emph{Combinatorics of
  representations of ${U}_q(\widehat{{sl}}(n))$ at $q=0$}, Comm. Math. Phys.
  \textbf{136} (1991), 543--566.

\bibitem{Kang03}
S.-J. Kang, \emph{Crystal bases for quantum affine algebras and combinatorics
  of {Y}oung walls}, Proc. London Math. Soc. (3) \textbf{86} (2003), no.~1,
  29--69.

\bibitem{KKM94}
S.-J. Kang, M.~Kashiwara, and K.~C. Misra, \emph{Crystal bases of {V}erma
  modules for quantum affine {L}ie algebras}, Compositio Math. \textbf{92}
  (1994), no.~3, 299--325.

\bibitem{KKMMNN91}
S.-J. Kang, M.~Kashiwara, K.~C. Misra, T.~Miwa, T.~Nakashima, and
  A.~Nakayashiki, \emph{Affine crystals and vertex models}, Infinite Analysis,
  Part A,B(Kyoto, 1991), 449--484, Adv. Ser. Math. Phys. {\bf16}, World Sci.
  Publ., River Edge, NJ (1992).

\bibitem{KKMMNN92}
\bysame, \emph{Perfect crystals of quantum affine {L}ie algebras}, Duke Math.
  J. \textbf{68} (1992), 499--607.

\bibitem{KL06}
S.-J. Kang and H.~Lee, \emph{Higher level affine crystals and {Y}oung walls},
  Algebr. Represent. Theory \textbf{9} (2006), no.~6, 593--632.

\bibitem{KP10}
S.-J. Kang and E.~Park, \emph{Quiver varieties and path realizations arising
  from adjoint crystals of type ${A}_n^{(1)}$}, Trans. Amer. Math. Soc.
  \textbf{363} (2011), no.~10, 5341--5366.

\bibitem{Kash90}
M.~Kashiwara, \emph{Crystalizing the $q$-analogue of universal enveloping
  algebras}, Comm. Math. Phys. \textbf{133} (1990), no.~2, 249--260.

\bibitem{Kash91}
\bysame, \emph{On crystal bases of the $q$-analogue of universal enveloping
  algebras}, Duke Math. J. \textbf{63} (1991), no.~2, 465--516.

\bibitem{KN94}
M.~Kashiwara and T.~Nakashima, \emph{Crystal graphs for representations of the
  $q$-analogue of classical {L}ie algebras}, J. Algebra \textbf{165} (1994),
  no.~2, 295--345.

\bibitem{KS97}
M.~Kashiwara and Y.~Saito, \emph{Geometric construction of crystal bases}, Duke
  Math. J. \textbf{89} (1997), no.~1, 9--36.

\bibitem{Kde09}
R.~Kodera, \emph{A generalization of adjoint crystals for the quantized affine
  algebras of type ${A}^{(1)}_n, {C}^{(1)}_n$ and ${D}^{(2)}_{n+1}$}, J.
  Algebraic Combin. \textbf{30} (2009), no.~4, 491--514.

\bibitem{L90}
G.~Lusztig, \emph{Canonical bases arising from quantized enveloping algebras},
  J. Amer. Math. Soc. \textbf{3} (1990), no.~2, 447--498.

\bibitem{L91}
\bysame, \emph{Quivers, perverse sheaves, and quantized enveloping algebras},
  J. Amer. Math. Soc. \textbf{4} (1991), no.~2, 365--421.

\bibitem{MM90}
K.~C. Misra and T.~Miwa, \emph{Crystal base for the basic representation of
  ${U}_q(sl(n))$}, Comm. Math. Phys. \textbf{134} (1990), no.~1, 79--88.

\bibitem{N94}
H.~Nakajima, \emph{Instantons on {A}{L}{E} spaces, quiver varieties, and
  {K}ac-{M}oody algebras}, Duke Math. J. \textbf{76} (1994), no.~2, 365--416.

\bibitem{N98}
\bysame, \emph{Quiver varieties and {K}ac-{M}oody algebras}, Duke Math. J.
  \textbf{91} (1998), no.~3, 515--560.

\bibitem{St02}
Y.~Saito, \emph{Crystal bases and quiver varieties}, Math. Ann. \textbf{324}
  (2002), no.~4, 675--688.

\bibitem{Sv06}
A.~Savage, \emph{Geometric and combinatorial realizations of crystal graphs},
  Algebr. Represent. Theory \textbf{9} (2006), no.~2, 161--199.

\bibitem{AS06}
A.~Schilling and P.~Sternberg, \emph{Finite-dimensional crystals ${B}^{2,s}$
  for quantum affine algebras of type ${D}^{(1)}_n$}, J. Algebraic Combin.
  \textbf{23} (2006), no.~4, 317--354.

\end{thebibliography}


\end{document}